\documentclass{article}

\usepackage{arxiv}

\usepackage[utf8]{inputenc}
\usepackage[T1]{fontenc} 
\usepackage{hyperref} 
\usepackage{url}
\usepackage{microtype}
\usepackage{doi}

\usepackage{graphicx}%
\usepackage{multirow}%
\usepackage{amsmath,amssymb,amsfonts}%
\usepackage{amsthm}%
\usepackage{mathrsfs}%
\usepackage[title]{appendix}%
\usepackage{textcomp}%
\usepackage{manyfoot}%
\usepackage{booktabs}
\usepackage{algorithm}%
\usepackage{algorithmicx}%
\usepackage{algpseudocode}%
\usepackage{listings}%
\usepackage[noabbrev,capitalize]{cleveref}
\usepackage{mathtools}
\usepackage{subcaption}
\usepackage{float}
\usepackage[table]{xcolor}
\usepackage{tikz}
\usetikzlibrary{positioning}
\usetikzlibrary{math}
\usepackage{cleveref}

\usetikzlibrary{fit,positioning}
\usetikzlibrary{backgrounds}
\usetikzlibrary{fit,backgrounds,calc}

\definecolor{light-gray}{RGB}{211,211,211}
\definecolor{lighter-gray}{RGB}{235,235,235}
\definecolor{lgreen}{RGB}{152,255,152}
\definecolor{lred}{RGB}{255,176,156}
\definecolor{blue1}{RGB}{144,202,249}
\definecolor{blue2}{RGB}{66,164,245}
\definecolor{blue3}{RGB}{30,120,229}

\renewcommand{\arraystretch}{1.3}
\newcolumntype{?}{!{\vrule width 1pt}}

\newtheorem{theorem}{Theorem}

\newtheorem{example}{Example}%
\newtheorem{remark}{Remark}%
\newtheorem{myalgorithm}[theorem]{Algorithm}

\newtheorem{definition}{Definition}%

\title{Solving Approximation Tasks \\ with Greedy Deep Kernel Methods}

\newif\ifuniqueAffiliation
\uniqueAffiliationfalse

\ifuniqueAffiliation 
\author{
	\And
	Coauthor \\
	\texttt{email} \\
	\AND
	Coauthor \\
	\texttt{email} \\
	\And
	Coauthor \\
	\texttt{email} \\
    \And
	Coauthor \\
	\texttt{email} \\
}
\else
\usepackage{authblk}

\newbox{\orcid}\sbox{\orcid}{\includegraphics[scale=0.06]{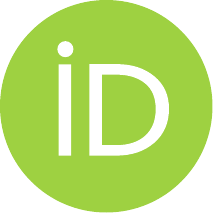}} 
\author[1]{%
	\href{https://orcid.org/0009-0002-8961-553X}{\usebox{\orcid}\hspace{1mm}Marian ~Klink\thanks{\texttt{marian.klink@ians.uni-stuttgart.de}}}%
}
\author[1]{%
	\href{https://orcid.org/0009-0000-4378-1757}{\usebox{\orcid}\hspace{1mm}Tobias ~Ehring\thanks{\texttt{tobias.ehring@ians.uni-stuttgart.de}}}%
}
\author[1]{%
	\href{https://orcid.org/0009-0009-1240-9229}{\usebox{\orcid}\hspace{1mm}Robin ~Herkert\thanks{\texttt{robin.herkert@ians.uni-stuttgart.de}}}%
}
\author[1]{%
	\href{https://orcid.org/0009-0004-9897-6614}{\usebox{\orcid}\hspace{1mm}Robin ~Lautenschlager\thanks{\texttt{robin.lautenschlager@ians.uni-stuttgart.de}}}%
}
\author[1]{%
	\href{https://orcid.org/0000-0002-1552-497X}{\usebox{\orcid}\hspace{1mm}Dominik ~Göddeke\thanks{\texttt{dominik.goeddeke@ians.uni-stuttgart.de}}}%
}
\author[1]{%
	{Bernard ~Haasdonk\thanks{\texttt{bernard.haasdonk@ians.uni-stuttgart.de}}}%
}
\affil[1]{Institute of Applied Analysis and Numerical Simulation and Stuttgart Center of Simulation Science (SC SimTech), University of Stuttgart, Stuttgart, Germany}
\fi

\renewcommand{\shorttitle}{Solving Approximation Tasks with Greedy Deep Kernel Methods}

\begin{document}
\maketitle

\begin{abstract}
Kernel methods are versatile tools for function approximation and surrogate modeling. In particular, greedy techniques offer computational efficiency and reliability through inherent sparsity and provable convergence.
Inspired by the success of deep neural networks and structured deep kernel networks, we consider deep, multilayer kernels for greedy approximation.
This multilayer structure, consisting of linear kernel layers and optimizable kernel activation function layers in an alternating fashion, increases the expressiveness of the kernels and thus of the resulting approximants. 
Compared to standard kernels, deep kernels are able to adapt kernel intrinsic shape parameters automatically, incorporate transformations of the input space and induce a data-dependent reproducing kernel Hilbert space. 
For this, deep kernels need to be pretrained using a specifically tailored optimization objective.
In this work, we not only introduce deep kernel greedy models, but also present numerical investigations and comparisons with neural networks, which clearly show the advantages in terms of approximation accuracies.
As applications we consider the approximation of model problems, the prediction of breakthrough curves for reactive flow through porous media and the approximation of solutions for parameterized ordinary differential equation systems.
\end{abstract}



\section{Introduction}

Kernel-based machine learning methods have demonstrated their powerful prediction capabilities in a variety of fields such as interpolation, regression or classification \cite{hofmann2008kernel,campbell2002kernel}. The fundamental basis of all kernel methods is the concept of positive definite kernels,
which implicitly compute the inner product in a high-dimensional reproducing kernel Hilbert
space (RKHS) between two, by means of a fixed non-linear feature map, transformed data points.
Based on this so-called kernel trick, kernel methods can be understood as linear algorithms on
transformed input data that are nevertheless able to model complex and nonlinear relationships,
making them widely and efficiently applicable. 
In addition, kernel-based methods are built on a solid theoretical framework, centered around representer theorems \cite{Steinwart2008,bohn2019representer,scholkopf2001generalized}, error analysis and convergence statements \cite{wendland2005scattered} or various oracle inequalities in the case of support vector machines \cite{Steinwart2008}. While kernel-based methods have proven to be very effective for small scattered datasets, they also suffer from some limitations:
First, kernel methods are usually bound to the use of a fixed feature map to transform the data into the high-dimensional feature space, such that the kernel matrix is the sole data representation. Choosing the best kernel function and tuning its parameters remain critical yet non-trivial tasks that majorly impact the approximation quality. Alternatively, constructing a data-dependent feature map which then defines the kernel function again requires a priori knowledge and limits the flexibility.
Second, when applied to large datasets, the computational costs can quickly exceed the range of feasibility due to the need to compute, store and possibly invert the increasingly ill-conditioned kernel matrix. 

Conversely, Neural Networks (NNs) exhibit the conceptional advantage that they do not rely on a fixed feature map, but instead automatically learn hierarchical features during training. 
They are also capable of dealing with arbitrary complex and high-dimensional vectorial data simply by increasing and tuning the architectures. 
However, NNs usually require large training datasets to achieve accurate predictions, often lack understanding and interpretability and demand excessive hyperparameter tuning.

In this work, we aim at mitigating the two major limitations of classical kernel methods, namely the fixed feature map and the computational costs, by combining kernel greedy methods with the expressive computational power and flexibility of NNs.
In this context, several researchers have already addressed the connection between kernel methods and NNs. 
For example, the NN Gaussian process \cite{lee2017deep,matthews2018gaussian} and the neural tangent kernel \cite{jacot2018neural,golikov2022neural} draw the connection between specific kernel methods and NNs in the infinite width limit, opening a promising alternative interpretation of NNs in the rich theoretical framework of kernel methods.
More recently, researchers have also investigated the construction of deep kernel architectures. 
Exemplarily, deep kernels using nonlinear input transformations have been introduced \cite{wilson2016deep} and a deep kernel representer theorem has been proven \cite{bohn2019representer}.
Further, our previous research efforts include the construction, analysis and application of structured deep kernel networks \cite{wenzel2026ana,wenzel2021struct} and 
the Vectorial Kernel Orthogonal Greedy Algorithm (VKOGA) \cite{wirtz2013vectorial,Santin2021VKOGA}.
In addition, we have successfully combined a $2$-layer deep kernel with VKOGA \cite{herkert2024greedy,wenzel2024application}, which has also been analyzed in \cite{wenzel2024data}.

In detail, greedy methods have two key advantages: they produce surrogate models with sparse expansions, leading to computationally efficient evaluations, and provide rigorous, sometimes optimal convergence guarantees. 
In particular, the associated error bounds contain a dimension-independent factor, which mitigates the curse of dimensionality.
However, classical results rely on the assumption that the target function lies in the RKHS generated by a fixed kernel — a rarely verifiable assumption in practice, for example whenever the target function is only implicitly given via provided datasets.
To relax this requirement, the deep kernel architecture adapts the kernel to the data and induces a target-dependent RKHS.
The architecture alternates between linear kernel layers and trainable kernel activation layers acting on the single-dimensions. 
This contrasts with standard neural networks, where the activation function is typically chosen a priori. Once the deep kernel parameters have been learned, a greedy algorithm selects representative greedy centers and the expansion coefficients follow from interpolation conditions.

In this work, we extend the work on $2$L-VKOGA \cite{wenzel2024application,wenzel2024data,herkert2024greedy} by introducing and examining a combination of VKOGA with even deeper kernels, in particular with up to $8$ layers.
Then, we systematically compare the resulting deep VKOGA models against NNs across a variety of different applications. These include model problem approximation tasks, the approximation of breakthrough curves of reactive flow through 3D porous geometry data as a continuation of \cite{herkert2024greedy}, and the approximation of solutions from parameterized ordinary differential equations (ODEs) such as the Lotka-Volterra ODE and Brusselator ODE.
In this series of numerical experiments, we demonstrate that the introduced deep VKOGA models frequently surpass the accuracy of ReLU NNs and graph NNs (GNNs) under similar and often smaller computational costs, underlining their versatility and reliability.
Further, our numerical experiments indicate that deeper kernel architectures are beneficial for complex approximation tasks.

The remainder of this work is organized as follows: In \cref{sec:2.1}, we elaborate on the construction of feedforward ReLU NNs. 
In \cref{sec:2.2}, we provide the theoretical background on greedy kernel interpolation and in particular the VKOGA method. Then, in \cref{sec:3}, we explain the extension of greedy kernel interpolation using deep kernels, including a brief discussion regarding training and analysis. In \cref{sec:5}, we explain the application and present the conducted numerical experiments for the three problem classes: Function approximation, breakthrough curve approximation and parameterized ODE solution approximation. Finally, we conclude this work in \cref{sec:6} with a brief summary and outlook to future extensions and improvements.

\section{Background}\label{sec:2}
In this section we provide the essential background on NNs and greedy kernel approximation on which the deep kernel greedy approach will build.

\subsection{Neural Networks}\label{sec:2.1}
Neural networks (NNs) represent classes of functions that are defined by a system of artificial neurons \cite{Sun2020,Murphy2022}. These function classes were originally inspired by the neural structure of the human brain. There are various types of NNs, depending on their neuron architectures, e.g., fully connected NNs, convolutional NNs and recurrent NNs \cite{goodfellow2016deep}. In this work, we consider feed forward fully connected NNs, where all neurons of adjacent layers are connected and the output of each layer is the input of the following layer. These types of NNs define a mapping from the input layer $0$ with $d_0 \in \mathbb{N}$ features to the output layer $L \in \mathbb{N}$ with $d_L \in \mathbb{N}$ features. The layers $\ell$ with $\ell \in \{1,\dots,L-1 \}$ and $d_\ell \in \mathbb{N}$ features are called hidden layers. NNs with more than one hidden layer are called deep neural networks. Two adjacent layers are recursively connected by
\begin{equation*}
   h_\ell:\mathbb{R}^{d_{\ell-1}}\longrightarrow\mathbb{R}^{d_\ell},
  \qquad h_\ell(x) ~=~ \sigma_\ell\left(b_\ell + W_\ell x \right)
\end{equation*}
for each layer index $\ell\in\{1,\dots,L\}$ \cite{Sun2020}. The trainable parameters $W_\ell \in \mathbb{R}^{d_{\ell}\times d_{\ell-1}}$ and $b_\ell \in \mathbb{R}^{d_\ell}$ are the weight matrix and the bias vector, respectively. Further, the function $\sigma_\ell:\mathbb{R}^{d_{\ell}}\longrightarrow\mathbb{R}^{d_\ell}$ is a component-wise activation function, e.g. rectified linear unit (ReLU) \cite{nair2010rectified}. For the output layer, that is for $\ell = L$, we use the identity activation $\sigma_L(x) \coloneqq x$ instead. 

When we assume the target function $f:\mathbb{R}^{d_0} \to \mathbb{R}^{d_L}$ to be a black-box, the NN surrogate model is given by
\begin{equation*}
    f(x) ~\approx~ h^{(L)}(x;\Theta) \coloneqq h_L \circ \dots \circ h_1 (x).
\end{equation*}
Here,
\[
\Theta \coloneqq (\text{vec}(W_1)^T,b_1^T, \text{vec}(W_2)^T,\dots,b_{L-1}^T,\text{vec}(W_{L})^T,b_L^T)^T \in \mathbb{R}^{N_p},\quad N_p =\sum_{\ell=1}^L \left(d_{\ell-1}d_\ell + d_\ell\right)
\]
is the parameter vector containing all optimizable NN parameters, where $\text{vec}(\cdot)$ is the vectorization operator flattening a matrix into a single column-vector.
In general, training the NN is purely data-driven and the goal is to find an optimal parameter vector $\Theta^\ast \in \mathbb{R}^{N_p}$ by solving an optimization problem, such as minimizing
\begin{eqnarray*}
    \mathcal{L}_{\mathrm{data}}(X_{\mathrm{train}},f(X_{\mathrm{train}}),\Theta) ~&\coloneqq&~ \frac{1}{N_{\mathrm{train}}}\sum_{j=1}^{N_{\mathrm{train}}}||f(x_j) ~-~ h^{(L)}(x_j;\Theta)||^2_2\\
    \Theta^\ast  ~&\coloneqq&~ \arg\min\limits_{\Theta \in \mathbb{R}^{N_p}} \mathcal{L}_{\mathrm{data}}(X_{\mathrm{train}},f(X_{\mathrm{train}}),\Theta)
\end{eqnarray*}
where $\mathcal{L}_{\text{data}}$ is the loss function on the training data $X_{\mathrm{train}} \coloneqq \{x_1,\dots,x_{N_{\mathrm{train}}} \} \subset \mathbb{R}^{d_0}$ with $N_{\mathrm{train}}$ the number of training points \cite{Sun2020}. 

\subsection{Greedy Kernel Approximation}\label{sec:2.2}
For a non-empty set $\Omega \subset \mathbb{R}^d$ and $m\in\mathbb{N}$,
a function
\[
  k:\Omega\times\Omega\longrightarrow\mathbb{R}^{m\times m}
\]
is a matrix-valued kernel if
$
  k(x,x')=k(x',x)^{\top}
$ for all $x,x'\in\Omega$ \cite{wittwar2018interpolation}.
For $X_N\coloneqq\{x_1,\dots,x_N\}\subset\Omega$ with $x_1,\dots,x_N$ pairwise distinct, we define the kernel matrix
\[
  K_{X_N}\;\coloneqq\;\bigl(k(x_i,x_j)\bigr)_{i,j=1}^{N}\in\mathbb{R}^{Nm\times Nm},
\]
which is to be understood as a block-matrix.
The kernel~$k$ is called positive definite (p.d.) if $K_{X_N}$ is positive semi-definite for every such~$X_N$ and
strictly positive definite (s.p.d.) if $K_{X_N}$ is positive definite for every such~$X_N$.
By Theorem 2.4 in \cite{wittwar2018interpolation} every p.d.\ kernel $k$ induces a unique RKHS
$  \mathcal{H}_k\;\coloneqq\;\mathcal{H}_k(\Omega)
  \subset\bigl\{f:\Omega\to\mathbb{R}^m\bigr\}$ 
such that $k(x,\cdot)a\in\mathcal{H}_k$ for all $x\in\Omega, a \in \mathbb{R}^m$, and
 the reproducing property
        $\langle f,k(x,\cdot)a\rangle_{\mathcal{H}_k}=f(x)^\top a$ holds for every $f\in\mathcal{H}_k$ and every $x\in\Omega,a \in \mathbb{R}^m$. In addition, a kernel is called translation-invariant if $k(x,x')=k(x+z,x'+z)$ for all $x,x'\in\Omega$ and $z \in \mathbb{R}^d$ such that $x+z, x'+z \in \Omega$.
On $\Omega\subset\mathbb{R}^{d}$, common choices for strictly positive definite and translation-invariant kernels are defined by radial basis functions (RBFs) depending on a shape parameter $\varepsilon>0$. Frequently used, shallow, scalar-valued kernel examples are given in Table \ref{tab:RBF_kernels}, which can easily be extended to matrix-valued kernels by $k(x,x'):=I_m  \cdot k_\text{RBF}(x,x')$, where $I_m$ denotes the $m \times m$ identity matrix.

\begin{table}[ht]
\centering
\renewcommand{\arraystretch}{1.2}
\begin{tabular}{|c|c|}
\hline
\textbf{Kernel} & \textbf{Closed form $k_\text{RBF}(x,x';\varepsilon)$}\\\hline
Gaussian &
$\displaystyle \exp\!\bigl(-\varepsilon^{2}\|x-x'\|_{2}^{2}\bigr)$ \\[2mm]
linear Matérn &
$\displaystyle\Bigl(1+\varepsilon\|x-x'\|_{2}\Bigr)\exp\!\bigl(-\varepsilon\|x-x'\|_{2}\bigr)$ \\[2mm]
quadratic Matérn &
$\displaystyle\Bigl(1+\varepsilon\|x-x'\|_{2}+\tfrac13\varepsilon^{2}\|x-x'\|_{2}^{2}\Bigr)
             \exp\!\bigl(-\varepsilon\|x-x'\|_{2}\bigr)$\\\hline
\end{tabular}
\caption{Selected strictly positive definite, translation-invariant scalar-valued RBF kernels.}
\label{tab:RBF_kernels}
\end{table}

\noindent Given a training set $X_{N_{\mathrm{train}}}\subset\Omega$ of pairwise distinct points and data
$  f(x_j) $
  for $j=1,\dots,N_{\mathrm{train}}$,
e.g., obtained from a function $f\in\mathcal{H}_k(\Omega)$, we seek an approximant
$s_{f,X_{N_{\mathrm{train}}}}\in\mathcal{H}_k(\Omega)$ that solves
the regularized least-squares regression problem
\begin{align}\label{eq:min_norm_problem}
  s_{f,X_{N_{\mathrm{train}}}} = \text{argmin}_{s\in\mathcal{H}_k(\Omega)}
  \left\{\sum_{j=1}^{N_{\mathrm{train}}} \Vert s(x_j)-f(x_j)\Vert^2_2 + \gamma \|s\|_{\mathcal{H}_k(\Omega)}^2 \right\} 
\end{align}
with hyperparameter $\gamma\geq 0$. 
Thanks to representer theorems a solution of problem~\eqref{eq:min_norm_problem} admits the finite representation
\begin{align}\label{eq:surrogateForm}
  s_{f,X_{N_{\mathrm{train}}}}(x)
  =\sum_{i=1}^{N_{\mathrm{train}}}\,k(x_i,x)\alpha_i,
\end{align}
where the coefficients $\alpha=(\alpha_1^\top,\dots,\alpha_{N_{\mathrm{train}}}^\top)^{\top} \in \mathbb{R}^{mN_{\mathrm{train}}}$ solve
$$
 (K_{X_{N_{\mathrm{train}}}} + \gamma I_{N_{\mathrm{train}}\cdot m} )\alpha=y
$$
with $y = (f(x_1)^\top,f(x_2)^\top,\ldots,f(x_{N_{\mathrm{train}}} )^\top)^\top \in \mathbb{R}^{m N_{\mathrm{train}}}$.
If $k$ is s.p.d.\ the Gramian matrix $K_{X_{N_{\mathrm{train}}}}$ is non-singular and if additionally $\gamma=0$,  we obtain a unique interpolant. For the scalar case this is proven in \cite[Chap.~16]{wendland2005scattered} and directly transfers to our diagonal matrix-valued kernel case. 

\noindent For large $N_{\mathrm{train}}$, solving the dense system and evaluating~\eqref{eq:surrogateForm} becomes very expensive and may suffer from ill-conditioning.  
Therefore, one usually selects a much smaller subset
$X_{n_\text{max}}\subset X_{\mathrm{train}}$ with $n_\text{max}\ll N_{\mathrm{train}}$.  
Greedy algorithms -- such as VKOGA~\cite{Santin2021VKOGA,wirtz2013vectorial} -- expand $X_n$ iteratively according to prescribed selection rules.  
A popular choice is the $f$-greedy rule, which at iteration $n\ge0$ selects the next greedy center
\begin{equation}\label{eq:f_greedy}
  x_{n+1}
  :=\underset{x\in X_{N_\mathrm{train}}}{\arg\max}
  \bigl\|f(x)-s_{f,X_n}(x)\bigr\|_2,\qquad
  X_{n+1}:=X_n \cup\{x_{n+1}\},\quad
  X_0:=\emptyset.
\end{equation}
Exemplarily, \cref{alg:VKOGA} visualizes the pseudocode of f-greedy VKOGA.

\begin{myalgorithm}[f-greedy VKOGA]\label{alg:VKOGA}
    \phantom{.}
    \hrule
    \vspace{1mm}
    \hrule
    \noindent\fcolorbox{white}{white}{
    \begin{minipage}{\dimexpr\linewidth-2\fboxsep-2\fboxrule}
    \begin{algorithmic}[1]
    \State \textbf{Input:} $k: \Omega \times \Omega \to \mathbb{R}^{m \times m}$ (strictly) positive definite kernel, 
    $f \in \mathcal{H}_k(\Omega)$ target function, \\
    $X_{N_\mathrm{train}} \subset \Omega$ training data, 
    $n_\text{max} \leq N_\mathrm{train}$ maximum number of greedy iterations, \\
    tolerance $\text{tol}$, selection rule $\eta(x; s) \coloneqq \Vert f(x) - s(x) \Vert_2$
    \vspace{2mm}
    \State \textbf{Initialize:} $X_0 \coloneqq \emptyset$, $s_{f,X_0} \coloneqq 0$, 
    $n = 0$
    \end{algorithmic}
    \end{minipage}}
    \hrule
    \noindent\fcolorbox{white}{lightgray!30}{
    \begin{minipage}{\dimexpr\linewidth-2\fboxsep-2\fboxrule}
    \begin{algorithmic}[1]
    \While{$n<n_\text{max}$ and $\underset{x\in X_{N_\mathrm{train}}}{\max} \eta\left(x, s_{f,X_n}\right) > \text{tol}$}
    \State Select $x_{n+1} :=\underset{x\in X_{N_\mathrm{train}}}{\arg\max} \eta\left(x, s_{f,X_n}\right)$
    \State Extend $X_{n+1}:=X_n \cup\{x_{n+1}\}$
    \State Compute $s_{f,X_{n+1}}$
    \State Set $n = n+1$
    \EndWhile
    \end{algorithmic}
    \end{minipage}}
    \hrule 
\end{myalgorithm}
\vspace{0.5cm}

Other strategies that result in sparse models are, e.g., $l_1$- penalization \cite{gao2010sparse} or support vector regression (SVR) \cite{scholkopf2002learning}.
Greedy schemes produce sparse surrogates and come with rigorous, sometimes even optimal, convergence guarantees \cite{santin2018greedy,santin2024optimality,wenzel2023analysis}. In particular, in the case of an s.p.d kernel and with the choice $\gamma = 0$, the associated error bounds contain a dimension-independent factor, thus mitigating the curse of dimensionality.  
For instance, for the scalar Gaussian kernel and the $f$-greedy procedure and $f \in \mathcal{H}_k(\Omega)$ holds
\[
  \min_{i=n+1,\dots,2n}
  \bigl\|f-s_{f,X_i}\bigr\|_{L^{\infty}(\Omega)}
  \;\le\;
  C\,e^{-c_1 n^{1/d}}\,
  \|f\|_{\mathcal{H}_k(\Omega)}\,n^{-1/2},
\]
which splits into a dimension-independent part $n^{-1/2}$ and a dimension-dependent part $e^{-c_1\,n^{1/d}}$.  
An analogous result holds for finitely smooth kernels \cite{wenzel2023analysis}.
\section{Greedy Deep Kernel Approximation}\label{sec:3}
In kernel-based approximation schemes, the choice of the
kernel is decisive.  
It determines (i) the associated RKHS, hence the class of
functions in which we search for a minimizer of \eqref{eq:min_norm_problem}, and (ii) the attainable approximation quality.
In practice, however, often neither the underlying RKHS nor the exact form of the target function are known
because the latter is often available only through data samples.
Consequently, selecting an optimal kernel is highly nontrivial.

\noindent Deep (or multi-layer) kernels mitigate this issue by learning suitable kernels from data:

\begin{enumerate}
  \item They can automatically adapt the shape parameters of basic kernels (e.g., RBF kernels) during training.
  \item They can incorporate affine transformations of the input domain.
  \item They can encode a mixture of several kernels and are therefore able to approximate composite targets whose components possess heterogeneous regularity.
\end{enumerate}
Overall, these design components raise the expressiveness of the corresponding kernel models and hence increase the chance that the unknown target function is contained or can be well-approximated in the RKHS of the corresponding deep kernel.
The final greedy deep kernel approximant then utilizes the deep kernel inside the greedy algorithm. For this, the deep kernel, which is defined by inner centers and trainable parameters, has to be trained over several epochs using a stochastic optimization algorithm. The complete pipeline is sketched in \cref{fig:DeepVKOGAFlowchart}.

\begin{figure}[H]
    \centering
\begin{tikzpicture}[every text node part/.style={align=center}]

        \node[rounded corners, fill=gray!30, minimum width=2cm, minimum height=0.75cm, align=center] (A) at (0,0) {\small Deep Kernel};
        \node[rounded corners, fill=gray!30, minimum width=2cm, minimum height=0.75cm, align=center] (B) at (4.0,0.0) {\small Stochastic Training};
        \node[rounded corners, fill=gray!30, minimum width=2cm, minimum height=0.75cm, align=center] (C) at (8.0,0.0) {\small Greedy Algorithm};

        \node (A1) at (0.0,1.5) {\small inner centers \\ \small $Z \subset X_{N_\text{train}}$};
        \node (A2) at (0.0,-1.5) {\small trainable parameters \\ \small $\Theta$};

        \node (B1) at (4.0,0.6) {\small for $i=1,\ldots,\text{epochs}$};
        \node (B2) at (4.0,-1.3) {\small update $\Theta$};

        \node (C1) at (8.0,0.6) {\small for $n=1,\ldots,n_\text{max}$};
        \node (C2) at (8.0,-1.6) {\small 1. Choose greedy center \\ \small $x_{n+1} \in X_{N_\text{train}}$ \\ \small 2. Update approximant};
        \node (Top1) at (0.0,2.0) {};
        \node (Bottom1) at (0.0,-2.2) {};
        \node (Top2) at (8.0,2.0) {};
        \node (Bottom2) at (8.0,-2.2) {};

    
    

        \begin{scope}[on background layer]
            \node[
                draw=gray!100,
                fill=gray!5,
                rounded corners,
                inner sep=5pt,
                fit=(A) (B) (A1) (A2) (B1) (B2)(Top1)(Bottom1)
            ] (G1) {};

            \node[anchor=south west] at ([xshift=-6pt,yshift=2pt]G1.north west)
                {\small\bfseries Deep Kernel Training};

            \node[
                draw=gray!100,
                fill=gray!5,
                rounded corners,
                inner sep=5pt,
                fit=(C) (C1) (C2)(Top2)(Bottom2)
            ] (G2) {};

            \node[anchor=south west] at ([xshift=-6pt,yshift=2pt]G2.north west)
                {\small\bfseries Approximant Generation};
        \end{scope}

        \draw[->] (A) -- (B);
        \draw[dashed, ->] (A1) -- (A);
        \draw[dashed, ->] (A2) -- (A);
        \draw[->] (B) -- (C);
        \draw[->] (5.5,-0.4) to[out=-40,in=-140] (2.5,-0.4);
        \draw[->] (9.5,-0.4) to[out=-40,in=-140] (6.5,-0.4);

        \node (Z) at (4.0,-2.5) {};
    \end{tikzpicture}
    \caption{Sketch of the greedy deep kernel approximation procedure: First, the inner centers and the trainable parameters are initialized. Second, the trainable parameters of the deep kernel are optimized by performing several stochastic batch optimization epochs. Third and final, the trained deep kernel is used inside the greedy algorithm (VKOGA), where in each greedy iteration a new greedy center is chosen and the final approximant is updated.}
    \label{fig:DeepVKOGAFlowchart}
\end{figure}
\noindent In the following, we first introduce the concept of a deep kernel consisting of linear and nonlinear kernel layers before we provide details on the choice of the inner centers and the training of the trainable parameters.

We formalize an $L$-layer deep kernel as a composition of
kernel-induced feature maps.

\begin{definition}[Deep kernel]\label{def:deep-kernel}
Let $\Omega\subseteq\mathbb{R}^{d_0}, L\in \mathbb{N}, L \geq 2$ and $M\in \mathbb{N}$ be fixed.
For each layer index $\ell\in\{1,\dots,L\}$
let
\[
  k_\ell:\mathbb{R}^{d_{\ell-1}}\times\mathbb{R}^{d_{\ell-1}}
  \;\longrightarrow\;\mathbb{R}^{d_\ell\times d_\ell}
\]
be a p.d. matrix-valued kernel.
For $\ell\in\{1,\dots,L-1\}$ choose inner centers $\{z_i^{\,\ell}\}_{i=1}^{M}\subset\mathbb{R}^{d_{\ell-1}}$ and define the layer map
\[
  f_\ell:\mathbb{R}^{d_{\ell-1}}\longrightarrow\mathbb{R}^{d_\ell},
  \qquad
  f_\ell(\cdot)
  \;:=\;
  \sum_{i=1}^{M}k_\ell\!\bigl(z_i^{\,\ell},\cdot\bigr)\,\alpha_{\ell,i},
  \qquad
  \alpha_{\ell,i}\in\mathbb{R}^{d_\ell}\;\text{(trainable parameters)}.
\]
With the propagating function
$
  F_{L-1}\ \coloneqq f_{L}\circ f_{L-1}\circ\cdots\circ f_{1}
$
we obtain the $L$-layer deep kernel
\begin{equation}\label{eq:deep-kernel}
  k^{(L)}:\Omega\times\Omega\longrightarrow\mathbb{R}^{d_L\times d_L},
  \qquad
  k^{(L)}(x,x')
  \;:=\;
  k_L\!\bigl(F_{L-1}(x),\,F_{L-1}(x')\bigr).
\end{equation}
\end{definition}

\begin{remark}
Because the composition of a p.d.\ kernel with a feature map
preserves positive definiteness, the kernel $k^{(L)}$ in
\eqref{eq:deep-kernel} is itself p.d.\ and, therefore, generates a unique
RKHS.  Consequently, $k^{(L)}$ can be deployed in the greedy center
selection strategy~\eqref{eq:f_greedy} to produce an approximant by selecting suitable greedy centers
$X_M := \{x_i\}_{i=1}^{M}\subset X_{\mathrm{train}}$ and linearly combining the deep kernel translates as in \cref{eq:surrogateForm}. Thus, the deep kernel approximant can be written as 
\begin{align}\label{eq:deepkernelapproximant}
  s_{f,X_{M}}(x)
  =\sum_{i=1}^{M}\,k^{(L)}(x_i,x)\alpha_i,
\end{align}
\end{remark}

Definition~\ref{def:deep-kernel} leaves three design variables unspecified:
(i)~the layer kernels $k_\ell, \ell = 1,\ldots,L$,
(ii)~the inner centers $\{z_i^{\,\ell}\}, 1 \leq \ell \leq L-1$ and
(iii)~the coefficient matrices $A_\ell=[\alpha_{\ell,1},\dots,\alpha_{\ell,M}] \in \mathbb{R}^{d_\ell\times M}$.
In the following we will comment on how we choose them. 

First, building on the analogy to feed-forward neural networks (see \cite{wenzel2023deep, wenzel2026ana} for the original
derivation and numerical results), 
we employ an architecture that {alternates two complementary
kernel layers. All odd-indexed layers are realized as \textit{linear kernel} layers; they enforce purely linear relationships between the components via a matrix-valued kernel.
Conversely, each even-indexed layer employs a \textit{componentwise} (block-diagonal) kernel activation, thereby introducing the non-linearities that are indispensable for expressive, universal function approximation. Moreover, it is sufficient to restrict $L$ to an even number, since for odd $L$ the final linear layer would not increase the expressiveness of the deep kernel as the RKHS would remain unchanged.

\paragraph{Linear kernel layer.}
Let $\ell\in\{1,3,5,\dots, L-1\}$ be an odd-layer index.
We employ the matrix-valued linear kernel
\begin{equation}\label{eq:linear kernel}
  k^{\mathrm{lin}}_{\ell}(z,z')
  \;:=\;
  \langle z,z'\rangle_{\mathbb{R}^{d_{\ell-1}}}\,
  I_{d_\ell},
  \qquad
  z,z'\in\mathbb{R}^{d_{\ell-1}}.
\end{equation}
Kernel~\eqref{eq:linear kernel} is  p.d. (but not s.p.d.) and induces the layer map
\[
  f_\ell:\mathbb{R}^{d_{\ell-1}}\longrightarrow\mathbb{R}^{d_\ell},
  \qquad
  f_\ell(\cdot)
  =\sum_{i=1}^{M}
     k^{\mathrm{lin}}_{\ell}\!\bigl(z_i^{\,\ell},\cdot\bigr)\,
     \alpha_{\ell,i}.
\]
As observed in \cite{wenzel2026ana}, the expansion collapses to a bias-free linear transformation,
\begin{equation}\label{eq:linear-layer}
  f_\ell(z)=W_\ell z,
  \qquad
  W_\ell\in\mathbb{R}^{d_\ell\times d_{\ell-1}},
\end{equation}
if and only if the span of the center points $\{z_i^\ell\}_{i=1}^M  \subset \mathbb{R}^{d_{l-1}}$ is a superset of the row
space of the matrix $W_\ell$. 
Additive biases are not required if we employ translation-invariant kernels
(e.g., Gaussian or Matérn) in the subsequent layer; such kernels depend only on differences of their arguments, rendering constant shifts
irrelevant to the layer output. Note that for training, we use representation (\ref{eq:linear-layer}) of the linear layer and consider the weight matrices $W_\ell$ as trainable parameters for the linear kernel layers. This means that the trainable parameters of the layer are exactly the entries of the weight matrix~$W_\ell$. 

\paragraph{Kernel activation layer.}
For every even layer index $\ell\leq L$ we introduce non–linearity
through a matrix-valued kernel that acts on the single dimensions
\begin{align}\label{eq:single-dim-kernel}
  &k^\text{act}_{\ell}:
    \mathbb{R}^{d_{\ell-1}}\times\mathbb{R}^{d_{\ell-1}}
    \;\longrightarrow\;
    \mathbb{R}^{d_\ell\times d_\ell} \textrm{ with } d_\ell=d_{\ell-1} \textrm{ by
construction}\\
  &k^\text{act}_{\ell}(z,z')
  \;:=\;
  \operatorname{diag}\!\bigl(
      k^\text{s}_\ell\bigl(z^{(1)},z'^{(1)}\bigr),
      \dots,
      k^\text{s}_\ell\bigl(z^{(d_\ell)},z'^{(d_\ell)}\bigr)
    \bigr) \qquad \text{for } z,z'\in \mathbb{R}^{d_{\ell-1}},
\end{align}
where $z^{(j)}$ denotes the $j$-th component of $z$.  In principle, one could assign a
distinct scalar kernel to each diagonal entry, thereby treating the
coordinate directions heterogeneously; this extension is not pursued
here.

\begin{remark}
Even if the scalar kernel $k^\text{s}_\ell$ is strictly positive definite, the
kernel \eqref{eq:single-dim-kernel} is in general only positive
definite: for vectors $z\neq z'$ it is possible that
$z^{(j)} = z'^{(j)}$ for some $j$, leading to zero eigenvalues of the
block-diagonal kernel matrix \cite{wenzel2026ana}.
\end{remark}
\noindent The associated layer map reads
\[
  f_\ell:\mathbb{R}^{d_{\ell-1}}\longrightarrow\mathbb{R}^{d_\ell},
  \qquad 
  f_\ell(z)
  =\sum_{i=1}^{M}
     k^\text{act}_{\ell}\!\bigl(z,z_i^{\,\ell}\bigr)\,
     \alpha_{\ell,i},
\]
for even $\ell \leq L-2$ with trainable column vectors
$\alpha_{\ell,i}\in\mathbb{R}^{d_\ell}$.  Collecting these columns
\[
  A_\ell
  \coloneqq
  \bigl[\alpha_{\ell,1}\; \alpha_{\ell,2}\;\dots\;\alpha_{\ell,M}\bigr]
  \in\mathbb{R}^{d_\ell\times M},
\]
layer $\ell$ contributes $M d_\ell$ free parameters to the model.

The propagating function corresponding to a deep kernel with linear kernel layers and kernel activation layers is sketched in \cref{sketch:Propagating_function}.
\begin{figure}[ht]
\centering
\begin{tikzpicture}[
  input/.style={circle, draw=green!70!black, fill=green!20, minimum size=0.5cm},
  hidden/.style={circle, draw=blue!70!black, fill=blue!20, minimum size=0.5cm},
  output/.style={circle, draw=red!70!black, fill=red!20, minimum size=0.5cm},
  every node/.style={font=\small}
]

\tikzmath{\x1 = 0; \x2 = 2.5; \x3 = 5; \x4 = 7.5; \x5 = 8.5; \x6 = 8; \x7 = 9.5; \x8 = 12;
          \y1 = -1; \y2 = -2; \ymiddle = -2.9; \y3 = -4; \y4 = -5;}

\node[input] (I1) at (\x1,\y2) {};
\node[draw=none] at (\x1,\ymiddle) {\scalebox{1.5}{$\vdots$}};
\node[input] (I2) at (\x1,\y3) {};

\node[hidden] (H11) at (\x2,\y1) {};
\node[hidden] (H12) at (\x2,\y2) {};
\node[draw=none] at (\x2,\ymiddle) {\scalebox{1.5}{$\vdots$}};
\node[hidden] (H13) at (\x2,\y3) {};
\node[hidden] (H14) at (\x2,\y4) {};

\node[hidden] (H21) at (\x3,\y1) {};
\node[hidden] (H22) at (\x3,\y2) {};
\node[draw=none] at (\x3,\ymiddle) {\scalebox{1.5}{$\vdots$}};
\node[hidden] (H23) at (\x3,\y3) {};
\node[hidden] (H24) at (\x3,\y4) {};

\node[hidden] (H31) at (\x4,\y1) {};
\node[hidden] (H32) at (\x4,\y2) {};
\node[draw=none] at (\x4,\ymiddle) {\scalebox{1.5}{$\vdots$}};
\node[hidden] (H33) at (\x4,\y3) {};
\node[hidden] (H34) at (\x4,\y4) {};

\node[draw=none] at (\x5,\ymiddle-0.15) {\scalebox{1.5}{$\ldots$}};


\node[hidden] (H51) at (\x7,\y1) {};
\node[hidden] (H52) at (\x7,\y2) {};
\node[draw=none] at (\x7,\ymiddle) {\scalebox{1.5}{$\vdots$}};
\node[hidden] (H53) at (\x7,\y3) {};
\node[hidden] (H54) at (\x7,\y4) {};

\node[output] (O1) at (\x8,\y2) {};
\node[draw=none] at (\x8,\ymiddle) {\scalebox{1.5}{$\vdots$}};
\node[output] (O2) at (\x8,\y3) {};

\node at (\x1 + 1.25, \y4 - 1.0) {$\underbrace{\hspace{2cm}}_{\parbox{3cm}{\centering Linear kernel \\ layer}}$};
\node at (\x2 + 1.25, \y4 - 1.0) {$\underbrace{\hspace{2cm}}_{\parbox{3cm}{\centering Kernel activation \\ layer}}$};
\node at (\x3 + 1.25, \y4 - 1.0) {$\underbrace{\hspace{2cm}}_{\parbox{3cm}{\centering Linear kernel \\ layer}}$};
\node at (\x7 + 1.25, \y4 - 1.0) {$\underbrace{\hspace{2cm}}_{\parbox{3cm}{\centering Linear kernel \\ layer}}$};

\node[align=center, above] at (\x1+0.8,-0.8) {$W_1 \in \mathbb{R}^{d_1 \times d_0}$};
\node[align=center, above] at (\x2+1.0,-0.4) {$A_2 \in \mathbb{R}^{d_2 \times M}$};
\node[align=center, above] at (\x3+1.0,-0.4) {$W_3 \in \mathbb{R}^{d_3 \times d_2}$};
\node[align=center, above] at (\x7+1.5,-0.8) {$W_{L-1} \in \mathbb{R}^{d_{L-1} \times d_{L-2}}$};

\foreach \i in {1,2}
  \foreach \j in {1,2,3,4}
    \draw[->] (I\i) -- (H1\j);

\foreach \i in {1,2,3,4}
    \draw[->] (H1\i) -- (H2\i);

\foreach \i in {1,2,3,4}
  \foreach \j in {1,2,3,4}
    \draw[->] (H2\i) -- (H3\j);

\foreach \i in {1,2,3,4}
    \draw[dashed,->] (H3\i) -- (\x4 + 0.75,\y\i);

\foreach \i in {1,2,3,4}
    \draw[dashed,->] (\x7 - 0.75,\y\i) -- (H5\i);


\foreach \i in {1,2,3,4}
  \foreach \j in {1,2}
    \draw[->] (H5\i) -- (O\j);
\vspace{1cm}
\end{tikzpicture}
\caption{Sketch of a propagating function $F_{L-1}: \mathbb{R}^{d_0} \to \mathbb{R}^{d_L-1}$ corresponding to an $L$-layer deep kernel $k^{(L)}:\mathbb{R}^{d_0} \times \mathbb{R}^{d_0}  \to \mathbb{R}^{d_L \times d_L}$. The input dimensions are marked in green, the hidden dimensions in blue and the output dimensions in red. The weight matrices $W_\ell$, $\ell$ odd, contain the optimizable parameters of the linear kernel layers. The weight matrices $A_\ell$, $\ell$ even, contain the optimizable parameters of the kernel activation layers.}
\label{sketch:Propagating_function}
\end{figure}

\paragraph{Inner centers.}
Next, we comment on the choice of the inner centers, i.e. those for $2\leq\ell \leq L-1$. Following the propagated‐center paradigm of
\cite{bohn2019representer,wenzel2026ana}, the inner centers are
not treated as independent variables.  
Instead, they are defined recursively by
\[
  z_i^{\,\ell}=F_{\ell-1}\bigl(z_i^{\,\ell-1}\bigr),
  \qquad
  \ell=2,\dots,L-1,\;\; i=1,\dots,M,
\]
with $F_{\ell-1}$ from~\eqref{eq:deep-kernel}, where the inner first-layer centers $\{z_i^1\}_{i = 1}^M$ are fixed before the training of the kernel layers. 
This implicit choice reduces the number of directly trainable
parameters while still allowing the inner centers to adapt by training the layer maps $f_\ell$.
Furthermore, this choice by propagation is optimal according to the deep kernel representer theorem \cite{bohn2019representer}.

\begin{example}
To illustrate the difference between deep and shallow kernels and the effect of the inner first-layer centers, \cref{fig:deep_kernel_visualization} visualizes different kernel functions based on the one-dimensional input domain $\Omega=[-3,3] \subset \mathbb{R}$. Precisely, the figure plots the kernel function $k^{(L)}(\cdot,1): [-3,3] \to \mathbb{R}_+$ where the second input of the kernel remains fixed. The figure visualizes the Gaussian kernel with $\epsilon=1$ and the linear Matérn kernel with $\epsilon=1$ from \cref{tab:RBF_kernels} and all $4$ different combinations of those two kernels that yield a $4$-layer deep kernel. For example, the $4$L Matérn-Gaussian kernel combines the linear Matérn kernel in the first kernel activation function layer with the Gaussian kernel in the second kernel activation function layer. The single dimensions are fixed to $d_\ell = 10$ for all linear and nonlinear layers and the respective weights are initialized randomly following the He initialization as described in \cite{he2015delving}. The inner first-layer centers vary from $M=1$ center at  $z_1^{1} = 0$ in the left plot, $M=1$ center at $z_1^{1} = 1$ in the middle plot and $M=3$ centers at $z_1^{1}=-1.5$, $z_2^{1}=0.0$ and $z_3^{1}=1.5$ in the right plot.
In all three plots, we see that the shallow Matérn and Gaussian kernels are axially symmetric at $x=1$ due to the second argument of the kernels being fixed at $1$. Contrarily, the deep kernels are axially symmetric at $x=0$ if the inner first-layer center is $z_1^{1} = 0$ (left plot) or at $x=1$ if the inner first-layer center is $z_1^{1} = 1$ (middle plot). When using multiple distinct inner first-layer centers, the deep kernels loose their axially symmetry completely as the nonlinear kernel activation function layer combines $M$ kernels linearly, which are all centered at different locations. Another interesting feature of the deep kernels is the vertical shift. While the shallow kernels converge to $0$ for $x \to \pm \infty$, the deep kernels are shifted vertically due to the bounded input transformation. In other words, the term $\Vert x - 1\Vert_2^2$ grows unbounded for $x \to \pm \infty$, whereas the term $\Vert F_L(x) - F_L(1) \Vert_2^2$ is bounded whenever $F_L$ is bounded. In our example, the propagating function $F_L$ is bounded due to the component-wise application of the nonlinear kernel activation function, which is a bounded function when using the linear Matérn and Gaussian kernels.
\begin{figure}
    \centering
    \input{deep_kernels_visualization.pgf}
    \caption{Visualization of shallow and $4$-layer deep kernels based on the Gaussian kernel with $\epsilon=1$ and the linear Matérn kernel with $\epsilon=1$ from \cref{tab:RBF_kernels}. Left: $M=1$ inner first-layer center $z_1^{1} = 0$. Middle: $M=1$ inner first-layer center $z_1^{1} = 1$. Right: $M=3$ inner first-layer centers $z_1^{1}=-1.5$, $z_2^{1}=0.0$ and $z_3^{1}=1.5$.}
    \label{fig:deep_kernel_visualization}
\end{figure}
\end{example}

\paragraph{Deep kernel training.}
Lastly, we explain how the weight matrices $W_\ell, \ell = 1,3,\ldots, L-1$ and coefficient matrices $A_\ell, \ell = 2,4,\ldots, L-2$ are optimized. 
For stochastic training of the deep kernel and hence the weight and coefficient matrices we split the training set  
\(
 X_{\mathrm{train}}
\)
into \(n_{\mathrm{batch}}\) disjoint batches
\(
  \mathcal{B}_p\subseteq X_{\text{train}}, p = 1,\ldots,n_{\mathrm{batch}}
\)
of sizes $|\mathcal{B}_p|=:B_p$.
On each batch we evaluate the leave-one-out (LOO) cross-validation
error of the current deep‐kernel approximant and minimize the aggregate
cost
\begin{equation}\label{eq:mini-batch-cost}
  C\!\bigl(\mathcal{B}_p,f(\mathcal{B}_p),\Theta\bigr)
  \;=\;
    \sum_{x_j\in\mathcal{B}_p}
  \bigl\|
    s_{f,\Theta, \mathcal{B}_p\setminus\{x_j\}}(x_j)-f(x_j)
  \bigr\|_{2}^{2},
\end{equation}
where \(s_{f,\Theta,\mathcal{B}_p\setminus\{x_j\}}\) denotes the
approximant built from the sub-batch
\(\mathcal{B}_p\setminus\{x_j\}\) by solving \eqref{eq:min_norm_problem} and 
\[
  \Theta\coloneqq (W_1,A_2, W_3,\dots,A_{L-2},W_{L-1})
\]
collects all trainable (coefficient) matrices. 
Gradient estimates of \eqref{eq:mini-batch-cost} with respect to
\(\Theta\) are obtained via automatic differentiation and used in a
batch gradient-descent scheme.

For an efficient implementation of the LOO-CV procedure on the current training batch $\mathcal{B}_p=\{x_1,\ldots,x_{B_p}\}$ we use an algorithm proposed by Rippa \cite{rippa1999algorithm}. 
We decompose $\mathcal{B}_p$ into $B_p$ sub-batches $\mathcal{B}_p \backslash \{x_j\}$ for $x_j \in \mathcal{B}_p$ and compute the kernel interpolant $s_{f,X_{{\mathcal{B}_p}\setminus \{x_j\}}}$ on $\mathcal{B}_p \backslash \{x_j\}$ 
as well as the error vector $e_j:=s_{f,X_{{\mathcal{B}_p}\setminus \{x_j\}}}(x_j) - f(x_j)$ for every $j=1,\ldots,B_p$.
Finally, we define the current cost as $C(\mathcal{B}_p, f(\mathcal{B}_p), \Theta) \coloneqq 1/B_p \sum_{j=1}^{B_p} \Vert e_j \Vert_2^2$. Fortunately, Rippa \cite{rippa1999algorithm} proves that the entries of the error vector can equivalently be computed by $e_j =  (K^{-1})_{jj} a_j$ with $a_j$ being the $j$-th block of $a$ where $a$ denotes the solution of the linear system $K a = y$ with $K$ being the kernel matrix of $\mathcal{B}_p$, $y$ being the corresponding target vector on the batch and 
$(K^{-1})_{jj}$ denoting the $j$-th diagonal block of $K^{-1}.$ Note that the blocks $(K^{-1})_{jj}$ can be computed from linear system solves without computing the inverse of the kernel matrix $K^{-1}$.
Similarly as explained in \Cref{sec:2.1}, in cases, where the kernel matrix $K$ is not positive definite, we consider the regularized kernel matrix $K_\gamma \coloneqq K + \gamma I $ using a regularization parameter $\gamma > 0$ instead. 
After the training of the kernel the global approximation task is addressed in the second stage, see the right hand side of \cref{fig:DeepVKOGAFlowchart}. This means combining the trained kernel with the VKOGA procedure \cite{Santin2021VKOGA,wirtz2013vectorial}, which results in the selection of the greedy centers $x_1, \ldots, x_{n_\text{max}}$ and the computation of the model's (outer) coefficient vectors as described in \Cref{sec:2.2}. At this greedy stage the trained kernel is treated as a fixed, shallow kernel.

\begin{remark}
    When $L=2$, Definition~\ref{def:deep-kernel} yields
\[
  k^{(2)}(x,x')
  =k_2\bigl(f_1(x),\,f_1(x')\bigr)
  =k_2\bigl(W_1x,\;W_1x'\bigr),
  \qquad
  W_1\in\mathbb{R}^{d_1\times d_0},
\]
so that the only trainable matrix \(W_1\) implements a global linear
transformation (including rotations, scalings, and shears) of the input
space \cite{wenzel2024data}.  This construction generalizes the
classical shape parameter of radial kernels to an arbitrary
affine deformation of the domain.

Further specialization \(d_0=d_1=1\) with
\(\displaystyle k_2(x,x')=\exp\bigl(-\varepsilon^2\|x-x'\|^2\bigr)\) (Gaussian kernel),
yields
\[
  k^{(2)}(x,x')
  =\exp\!\bigl(-\|W_1x - W_1x'\|^2\bigr)
  =\exp\!\bigl(-\varepsilon^2\|x-x'\|^2\bigr),
  \quad \text{with } W_1 = \varepsilon I,
\]
so that two layer-VKOGA (2L-VKOGA) coincides with greedy Gaussian kernel interpolation in
which the shape parameter \(\varepsilon\) is optimized in a
preprocessing step.
\end{remark}

\section{Numerical Experiments}\label{sec:5}

As applications for the deep VKOGA models and their comparison with NNs we examine three different problem classes. We conduct numerical experiments to approximate model problems in \cref{sec:5.1}, breakthrough curves of reactive flow through 3D porous geometries in \cref{sec:5.2} and solutions of parameterized ODEs in \cref{sec:5.3}.
The underlying Python code, which can be used to reproduce all numerical results, is provided via DaRUS\footnote{\label{note1}DaRUS dataset: \url{https://doi.org/10.18419/DARUS-5167}}. 
The VKOGA and two-layer VKOGA code is based on \cite{wenzel2024data}.
Further, we leverage the Python library PyTorch\footnote{PyTorch: \url{https://pytorch.org/}} to construct and train the NN and deep kernel architectures. 
In detail, we use the Adam optimizer\footnote{PyTorch Adam: \url{https://docs.pytorch.org/docs/stable/generated/torch.optim.Adam.html} \cite{kingma2014adam}} to tune all optimizable network parameters and perform a $5$ or $10$-fold cross-validation to tune problem-dependent hyperparameters such as learning rate, number of hidden neurons, number of layers, use of dropout in the case of NNs and the type of kernels in the case of deep VKOGA models. 
For simplicity, we always use the same single-dimensional kernel in all kernel activation layers of a deep kernel and the $f$-greedy selection criterion to select the greedy centers.
In all conducted experiments, we allow for early termination of the greedy iterations whenever the training residual falls below an absolute tolerance of $10^{-12}$. Further, we choose the inner first-layer centers of the deep kernel as a subset of the training data and propagate them successively through all deep kernel layers. 
The details of all cross-validations as well as the resulting model architectures are documented in the accompanying repository.

For all problem classes, we introduce the notation
\begin{align*}
    D_\text{train} &\subset \mathbb{R}^{d_\text{in}} \times \mathbb{R}^{d_\text{out}}, \quad | D_\text{train}| = N_\text{train}\\
    D_\text{test} &\subset \mathbb{R}^{d_\text{in}} \times \mathbb{R}^{d_\text{out}}, \quad | D_\text{test}| = N_\text{test}, \quad D_\text{test} \cap D_\text{train} = \emptyset
\end{align*}
for the train and test datasets, respectively, where $N_\text{train} \in \mathbb{N}$ and $N_\text{test}\in \mathbb{N}$ denote the number of train and test data samples and $d_\text{in}\in \mathbb{N}$ and $d_\text{out}\in \mathbb{N}$ represent the input and output dimension of the target functions. 
With this, we introduce the relative test error measure to judge the approximation capability of a surrogate model $s: \mathbb{R}^{d_\text{in}} \to \mathbb{R}^{d_\text{out}}$ via
\begin{align}\label{eq:e_rel}
    e^\text{rel} \coloneqq e^\text{rel}\left(D_\text{test}, s\right) \coloneqq \frac{1}{N_\text{test}} \sum_{(x,y) \in D_\text{test}} \frac{\Vert s(x) - y \Vert_2^2}{\Vert y \Vert_2^2} .
\end{align}
To measure the computational effort of all model architectures in the offline phase, that is the runtime during training, we average the runtime over $5$ independent offline runs. Similarly, we average the online runtime, that is the runtime of a single prediction, over all $N_\text{test}$ samples of the test dataset. 
All runtime measurements have been performed under similar conditions on a system using an AMD Ryzen 5 PRO 8540U CPU and 40 GB of RAM.

\subsection{Model Problems}\label{sec:5.1}

Inspired by \cite{wenzel2024data}, we start with three model problem approximation tasks. The target functions $f_2$, $f_3$ and $f_4$, where the subscripts are not an enumeration but indicate the input dimensions, are:
\begin{align*}
    f_2: \Omega \coloneqq [0,1]^{2} \to \mathbb{R}, \quad &f_2: x \mapsto \exp\left({-4\sum_{i=1}^{2}\left(x_i-0.5\right)^{2}}\right) \\
    f_3: \Omega \coloneqq [0,1]^{3} \to \mathbb{R}, \quad &f_3: x \mapsto \exp\left({-4\sum_{i=1}^{3}\left(x_i-0.5\right)^{2}}\right) + 2|x_1-0.5|\\
    f_4: \Omega \coloneqq [0,1]^{4} \to \mathbb{R}, \quad &f_4: x \mapsto \exp\left({-4\sum_{i=1}^{4}\left(x_i-0.5\right)^{2}}\right) + \exp\left({-9\sum_{i=1}^{2}\left(x_i-0.3\right)^{2}}\right)
\end{align*}
For all three target functions we choose $N_\text{train}=5000$ random uniformly distributed training samples and $N_\text{test}=1000$ random uniformly distributed test samples on $\Omega$ to create the train and test datasets.
Thus, the data availability per dimension decreases drastically while the target functions $f_2,f_3$ and $f_4$ inherit increasingly difficult characteristics, rendering the approximation task increasingly challenging.
 
For the deep VKOGA models, we choose $50$ inner first-layer centers and perform $50$ greedy iterations and $200$ training epochs with a batch size of $100$ and a regularization parameter $\gamma=0.001$ for the computation of Rippa's loss.
In the cross-validation, we choose between a linear Matérn kernel with shape parameter $\epsilon=1.0$ or $0.1$ or a Gaussian kernel with shape parameter $\epsilon=0.1$, $d_\ell = 10$ or $20$ single dimensions per layer $\ell=1,\ldots,L$ and a learning rate of $\lambda=0.01$ or $0.001$.

For the ReLU NNs, we perform $500$ training epochs with a batch size of $32$ or $100$ and use $10\%$ of the training data as validation data for early stopping. 
In the cross-validation, we choose between $10,20$ or $50$ hidden neurons per layer, a learning rate of $\lambda=0.01, 0.001$ or $0.0001$ and the use of dropout with dropout rate $0.2$ as boolean hyperparameter.

\begin{figure}[!htbp]
    \centering
    \input{MDtest_single_plot.pgf}
    \caption{Comparison between deep VKOGA models and NNs regarding the functions $f_2$ (left column), $f_3$ (middle column) and $f_4$ (right column). Top row: Mean relative test error evaluated on the fixed test dataset as a function of network depth (number of layers). The error bars indicate the $5$th-$95$th percentile range of relative test errors across the test dataset. The dashed lines correspond to the deep VKOGA models, the solid lines to the ReLU NNs. Middle row: Training loss, that is Rippa's loss in the case of deep VKOGA models and the MSE loss in the case of NNs, over the number of training epochs. Only every $10$-th loss evaluation is plotted to enhance visibility. Bottom row: Training residual norm over the number of selected greedy centers.}
    \label{fig:MD_test}
\end{figure}

The results of the conducted numerical experiments are displayed in \cref{fig:MD_test}, where the y-axis scales are not fixed.
The first row visualizes the relative test errors from \cref{eq:e_rel} with respect to different numbers of layers $L$, where the error-bars display the $10\%$ and $90\%$ percentile values of the relative test error. The second row visualizes the training loss with respect to the number of epochs and the third row the training residual $L_\infty$ norm of the VKOGA models with respect to the number of selected greedy centers. 
Here, $L=1$ for VKOGA is equivalent to the standard VKOGA model without the use of a deep kernel.
Moreover, we emphasize that, for example, the $4$L-VKOGA model contains the same number of linear layers as a $2$L-NN, namely $2$, but in addition an optimizable kernel activation function layer and the outer greedy kernel expansion layer. We only count layers containing optimizable parameters.

Examining the relative test error plots in \cref{fig:MD_test}, we clearly observe that the deep VKOGA models achieve significantly smaller relative test errors compared to the ReLU NNs, especially when considering similarly deep architectures. 
Whereas the $2$L-VKOGA model is the most accurate model for $f_2$, the deeper VKOGA models achieve smaller relative test errors for the target functions $f_3$ and $f_4$, demonstrating that deeper kernels can be beneficial for higher dimensional functions with different characteristics in each input dimension.
Note that the superior performance of the $2$L-VKOGA model for function $f_2$ is not surprising due to the fact that the target function is nothing else than the Gaussian kernel function with $\epsilon=2$ and centered at $x'=(0.5,0.5)^T$. Thus, a simple linear transformation is sufficient for the $2$L-VKOGA model using the Gaussian kernel with $\epsilon=0.1$ to achieve accurate predictions. This also emphasizes the importance of the shape parameter, since the $1$L-VKOGA model has not been able to reach similar accuracies with neither the Gaussian kernel with $\epsilon=0.1$ nor the Gaussian kernel with $\epsilon=1.0$.
Additionally, we can see that the most accurate VKOGA model is for all three target functions clearly more accurate than the most accurate ReLU NN. 
Taking a closer look at the corresponding training losses, \cref{fig:MD_test} (second row) illustrates Rippa's loss for the deep VKOGA models and the MSE loss for the ReLU NNs, evaluated at every $10$-th epoch. 
Here, we observe that in most cases a smaller training loss corresponds to a smaller relative test error. 
Whereas the $4$L-, $6$L- and $8$L-NNs achieve similar training losses after $500$ epochs, the deep kernels seem to profit from architectures with $4$ or more layers, at least regarding the more challenging target functions $f_3$ and $f_4$. 
However, we still observe a stagnation of the relative test error when increasing the number of layers beyond $4$, while the offline and online runtimes increase due to increasing model complexity, as discussed later with the results from \cref{tab:MD_test}.
Moreover, the deep VKOGA training losses often do not seem to stagnate after $200$ epochs, suggesting even better accuracies at the costs of more expensive offline phases.
Similarly, the $L_\infty$ norms of the deep VKOGA training residuals decrease nearly monotonically when increasing the number of selected greedy centers, as can be seen in the bottom row of \cref{fig:MD_test}. 
Again, we observe in most cases a direct correspondence between smaller training residual norms and smaller relative test errors. 
In addition, the results suggest that a larger number of greedy iterations can increase the accuracy of the deep VKOGA models even further, but increases simultaneously the offline and online costs.

We summarize the averaged runtime measurements for all models and all target functions in \cref{tab:MD_test}. 
As can be seen, the deep VKOGA models are in general slightly more efficient in the offline phase compared to the ReLU NNs using the same number of optimizable layers. The sole exceptions are $8$L-VKOGA for $f_2$ and $6$L-VKOGA and $8$L-VKOGA for $f_3$, where the offline phases are slightly more expensive compared to the corresponding NN architectures. 
Thus, we compare surrogate models requiring similar computational resources in the offline phase.
Regarding the online costs we often see slight advantages for the NNs, while the online costs of the deep VKOGA models are usually in similar orders of magnitude.

\begin{table}[!htbp]
    \centering
    \begin{tabular}{l?cc?cc?cc?}

        & \multicolumn{2}{c?}{\cellcolor{light-gray} $f_2$} & \multicolumn{2}{c?}{\cellcolor{light-gray} $f_3$} & \multicolumn{2}{c?}{\cellcolor{light-gray} $f_4$} \\

        & \cellcolor{lighter-gray} $t_\text{off}$ in [s] & \cellcolor{lighter-gray} $t_\text{on}$ in [s] & \cellcolor{lighter-gray} $t_\text{off}$ in [s] & \cellcolor{lighter-gray} $t_\text{on}$ in [s] & \cellcolor{lighter-gray} $t_\text{off}$ in [s] & \cellcolor{lighter-gray} $t_\text{on}$ in [s]\\
        \hline
        \hline
        \cellcolor{light-gray}1L-VKOGA & 1.00e-01 & 1.83e-06 & 1.03e-01 & 1.37e-06 & 1.02e-01 & 1.53e-06 \\ 

        \cellcolor{light-gray}2L-VKOGA & 6.90e+00 & 5.57e-04 & 8.22e+00 & 3.61e-06 & 7.54e+00 & 2.56e-06 \\ 

        \cellcolor{light-gray}4L-VKOGA & 2.42e+01 & 1.05e-05 & 2.06e+01 & 5.51e-06 & 2.22e+01 & 9.46e-06 \\ 

        \cellcolor{light-gray}6L-VKOGA & 2.88e+01 & 4.14e-06 & 4.39e+01 & 1.45e-05 & 4.68e+01 & 1.30e-05 \\ 

        \cellcolor{light-gray}8L-VKOGA & 6.32e+01 & 2.63e-05 & 5.90e+01 & 1.02e-05 & 6.64e+01 & 2.18e-05 \\
        \hline
        \hline
        \cellcolor{lighter-gray}2L-NN & 1.96e+01 & 2.89e-06 & 1.98e+01 & 7.73e-07 & 2.11e+01 & 9.18e-07 \\ 

        \cellcolor{lighter-gray}4L-NN & 2.56e+01 & 9.36e-07 & 3.62e+01 & 1.03e-06 & 2.56e+01 & 9.77e-07 \\ 

        \cellcolor{lighter-gray}6L-NN & 3.84e+01 & 1.01e-06 & 4.76e+01 & 1.13e-06 & 3.44e+01 & 9.92e-07 \\ 

        \cellcolor{lighter-gray}8L-NN & 4.99e+01 & 1.52e-06 & 6.02e+01 & 1.35e-06 & 5.56e+01 & 1.19e-06 \\ 
        \hline
        \hline
    \end{tabular}
	\caption{Averaged online and offline runtime measurements in seconds for all deep kernel models and NNs applied to functions $f_2$, $f_3$ and $f_4$.}
	\label{tab:MD_test}
\end{table}

Finally, we examine the effect of depth for fixed width in \cref{fig:MDtest_depth_vs_width}, that is we display the mean relative test error for deep VKOGA models and NNs over a different number of layers while we ensure that the number of trainable parameters remains nearly constant. 
For all models, we display the number of hidden neurons or single dimensions, respectively, and the total number of trainable parameters in \cref{tab:MD_test_width_vs_depth}. 
In the case of the deep VKOGA models, we use the Gaussian kernel with $\epsilon = 0.1$ for function $f_2$ and the linear Matérn kernel with $\epsilon=0.1$ for functions $f_3$ and $f_4$. Further, we always use a learning rate of $\lambda=0.01$, whereas all other parameters remain as specified around the discussion of \cref{fig:MD_test}.
In the case of the case of the NNs, we always use a learning rate of $\lambda=0.001$, a batch size of $32$ and no dropout, whereas again all other parameters remain as specified around the discussion of \cref{fig:MD_test}. In general, this choice of hyperparameters corresponds to the most common choice during the cross validations for the previous experiment. 
Examining \cref{fig:MDtest_depth_vs_width}, we again see the advantage of the examined deep VKOGA models compared to the NNs in terms of accuracy, with the sole exceptions being $2$L-VKOGA for $f_2$ and $f_3$. 
More interestingly, we see that network depth is really important for the deep kernels when approximating the more challenging functions $f_3$ and $f_4$. 
Here, the deeper kernels achieve in general more accurate predictions with nearly the same number of trainable parameters. 
In particular the difference in accuracy between the $2$L-VKOGA models, where no nonlinear kernel activation function is used, and the $4$L-VKOGA models is large, emphasizing the importance of nonlinear transformations for these two approximation tasks. 
Contrarily, for function $f_2$ we again see that the $2$L-VKOGA model is the most accurate model and that adding network depth renders the optimization problem during training only more complicated.
Examining the NNs, we see that an increased network depth does not yield more accurate predictions in these examples. For functions $f_2$ and $f_3$, the $2$L-NN is the most accurate NN model, whereas only for function $f_4$ the $6$L-NN yields slightly more accurate predictions.
\begin{figure}[!htbp]
    \centering
    \input{MDtest_width_vs_depth_error_rel_L.pgf}
    \caption{Mean relative test error of deep VKOGA models and NNs with a nearly constant number of trainable parameters evaluated on the fixed test dataset as a function of network depth (number of layers). The error bars indicate the $5$th-$95$th percentile range of relative test errors across the test dataset. The dashed lines correspond to the deep VKOGA models, the solid lines to the ReLU NNs.}
    \label{fig:MDtest_depth_vs_width}
\end{figure}

\begin{table}[!htbp]
\centering
    \begin{tabular}{l?cc?cc?cc?}

        & \multicolumn{2}{c?}{\cellcolor{light-gray} $f_2$} & \multicolumn{2}{c?}{\cellcolor{light-gray} $f_3$} & \multicolumn{2}{c?}{\cellcolor{light-gray} $f_4$} \\

        & \cellcolor{lighter-gray} $d_\ell$ & \cellcolor{lighter-gray} $N_p$ & \cellcolor{lighter-gray} $d_\ell$ & \cellcolor{lighter-gray} $N_p$ & \cellcolor{lighter-gray} $d_\ell$ & \cellcolor{lighter-gray} $N_p$\\
        \hline
        \hline
        \cellcolor{light-gray}2L-VKOGA & 910 & 1970 & 610 & 2030 & 460 & 2090 \\ 

        \cellcolor{light-gray}4L-VKOGA & 24 & 1974 & 24 & 2048 & 24 & 2122 \\ 

        \cellcolor{light-gray}6L-VKOGA & 14 & 1970 & 14 & 2034 & 14 & 2098 \\ 

        \cellcolor{light-gray}8L-VKOGA & 10 & 1970 & 10 & 2030 & 10 & 2090 \\ 

        \hline
        \hline
        \cellcolor{lighter-gray}2L-NN & 525 & 2101 & 420 & 2101 & 350 & 2101 \\ 

        \cellcolor{lighter-gray}4L-NN & 31 & 2109 & 31 & 2140 & 31 & 2171 \\ 

        \cellcolor{lighter-gray}6L-NN & 22 & 2113 & 22 & 2135 & 22 & 2157 \\ 

        \cellcolor{lighter-gray}8L-NN & 18 & 2125 & 18 & 2143 & 18 & 2161 \\ 
        \hline
        \hline
    \end{tabular}
	\caption{Number of hidden neurons or single dimensions $d_\ell$ per layer and number of total trainable parameters $N_p$ for every deep VKOGA model and NN displayed in \cref{fig:MDtest_depth_vs_width}.}
	\label{tab:MD_test_width_vs_depth}
\end{table}

\subsection{Breakthrough Curves}\label{sec:5.2}

\paragraph{The breakthrough curve approximation problem}
Next, we consider the approximation of breakthrough curves from chemical species flowing through porous geometries at the pore scale as a continuation of \cite{herkert2024greedy}. 
In detail, the breakthrough curves $a: \mathbb{R}^{+} \to \mathbb{R}^{+}$ are given by time-dependent integrals of the chemical species concentrations $c: \mathbb{R}^{3} \times \mathbb{R}_0^{+} \to \mathbb{R}_0^{+}$ over the outlet boundaries $\Gamma_\text{outlet} \subset \partial \Omega_\text{geo}$ of the 3D geometries $\Omega_\text{geo} \coloneqq [0,1]^{3}$, that is
\begin{align*}
    a(t) \coloneqq \int_{\Gamma_\text{outlet}} c(x,t) d\sigma, \quad t \in \mathbb{R}_0^{+}
\end{align*}
where $d\sigma$ represents a surface element of $\Gamma_\text{outlet}$. 
For this task, we consider an experimental data set with a total number of $N=59$ data samples. 
In detail, data sample $i$, $1 \leq i \leq 59$ consists of geometry data $x_i \in \{0,1\}^{2 \times 150 \times 150 \times 150}$ and the evaluation of the corresponding breakthrough curve $y_i \in \mathbb{R}^{500}$, $(y_i)_j \approx a(t_j)$ at $500$ prescribed timesteps $t_j$ for $j=1,\ldots,500$.
The geometry data are given by two boolean values at a uniform grid on $\Omega_\text{geo}$ with a total number of $N_v = 150^{3}$ voxel units, encoding the voxel to consist of either free pores, washcoat (porous catalyst) or solid substrate.
Exemplarily, \cref{fig:3Dgeoms} visualizes a 3D geometry on the left, a cross section in the middle and the resulting breakthrough curve on the right.

\begin{figure}[!htbp]
    \centering
    \includegraphics[width=\textwidth]{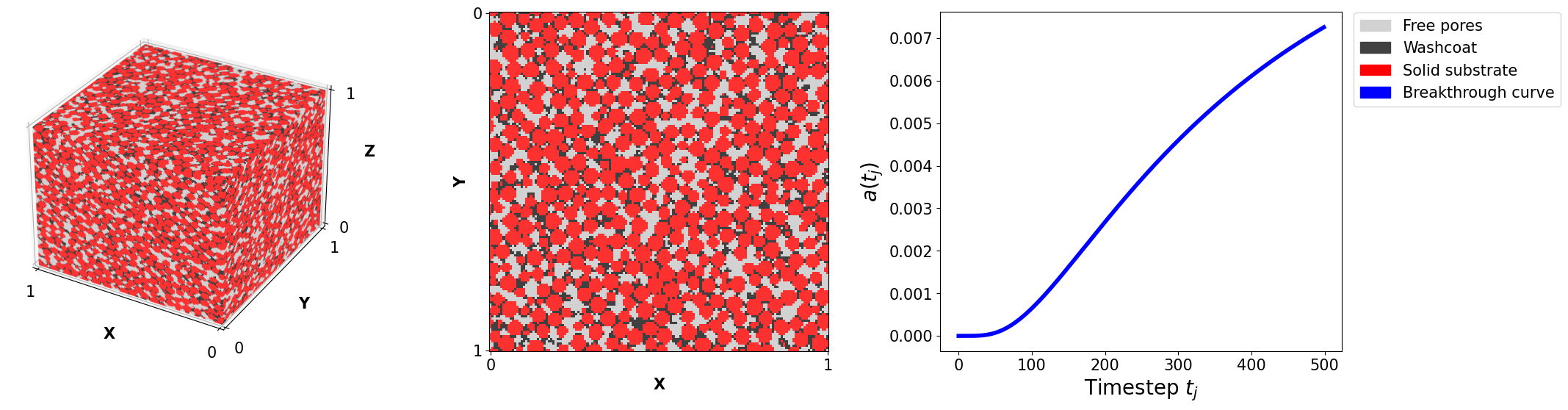}
    \caption{Visualization of a single data sample. Left: 3D geometry. Middle: Cross section of the $x$-$y$ plane at $z=0.5$. Right: Corresponding breakthrough curve.}
    \label{fig:3Dgeoms}
\end{figure}

Since dealing with the high-dimensional geometry data directly is computationally intensive, we make use of a principal component analysis (PCA) feature map to extract low-dimensional features a priori, as discussed in \cite{herkert2024greedy}.
For this, we first flatten the $N=59$ input tensors $x_i$ to obtain $N$ vectors $z_1,\ldots,z_N \in \mathbb{R}^{2 \cdot N_\nu}$, which define the matrix $Z \coloneqq [z_1,\ldots,z_N] \in \mathbb{R}^{2N_\nu \times N}$. 
Then, we define the PCA feature map $\Phi: \mathbb{R}^{2N_\nu} \to \mathbb{R}^{6}$ extracting $6$ features by inner products with the first $6$ left singular vectors of the singular value decomposition of $Z$
\begin{align*}
    Z = U \Sigma V^T, \quad U_6 = [u_1,u_2,\ldots,u_6] \in \mathbb{R}^{2N_\nu \times 6}, \quad \Phi(x) \coloneqq U_6^T x, \quad x \in \mathbb{R}^{2N_\nu}
\end{align*}
with $U = [u_1, \ldots, u_{2N_v}] \in \mathbb{R}^{2N_v \times 2N_v}$, $V \in \mathbb{R}^{N \times N}$ being the matrices with left and right singular vectors and 
$\Sigma \in \mathbb{R}^{2N_v \times N}$ being the matrix with singular values on its diagonal.
Consequently, we aim at approximating a target function $g: \mathbb{R}^{6} \times \mathbb{R} \to \mathbb{R}$, $(z,t) \mapsto a(t; z)$ after mapping the input tensors via the PCA feature map to the space $\mathbb{R}^{6}$.
For this, we distinguish between two different applications of the deep VKOGA models and NNs, which we refer to as discrete-time (DT) and continuous-time (CT) approach. 
In the case of the DT approach, we train a surrogate model to learn the mapping
\begin{align*}
    g: \mathbb{R}^6 \to \mathbb{R}^{500}, \quad z \mapsto [a(t_1;z),\ldots,a(t_{500};z)] .
\end{align*}
Thus, the surrogate model can predict the solution at all time grid points simultaneously. 
Alternatively, in the case of the CT approach, we train a surrogate model to learn the mapping
\begin{align}\label{eq:CT_surrogate}
    g: \mathbb{R}^6 \times \mathbb{R} \to \mathbb{R}, \quad (z,t) \mapsto a(t;z), 
\end{align}
that is the surrogate model predicts the solution at a single point in time. 
Here, the surrogate model is not restricted to the predefined time grid $\Omega_t^{\text{h}}$, both in terms of test and training datasets. 
However, due to the increase in the dataset sizes (by factor $500$ in this instance), the CT surrogate models are in general less efficient compared to the DT surrogate models in both the offline and online phase.

\paragraph{Numerical experiments}
Due to the small amount of available data, we randomly divide it into 3 train-test splits, where each split consists of $47$ training data samples and $12$ test data samples. 

For the DT deep VKOGA models, we choose $12$ inner first-layer centers and perform $47$ greedy iterations and $500$ training epochs with a batch size of $12$ and a regularization parameter $\gamma = 0.0001$ for the computation of Rippa's loss. Additionally, we perform a cross-validation to choose between the linear or quadratic Matérn kernel with shape parameter $\epsilon = 1.0$ or $0.1$, $d_\ell = 10$ or $20$ single dimensions per layer $\ell=1,\ldots,L$, a regularization parameter for the interpolation conditions of $0$ or $0.001$ and a learning rate of $\lambda=0.01$ or $\lambda=0.001$. 
Note, that performing $47$ greedy iterations corresponds to standard deep kernel interpolation on the complete training dataset, but with a more numerically stable algorithm. In fact, we allow the greedy algorithm to stop the greedy iterations early, whenever the norm of the training residual falls below a tolerance of $10^{-12}$. However, for this particular example of the breakthrough curve approximation, the stopping criterion is never reached and all $47$ available greedy centers are chosen.
Contrarily, due to the increased offline costs, we do not perform a cross validation for the CT deep VKOGA models but empirically choose a suitable architecture using the linear Matérn kernel with shape parameter $\epsilon=0.1$, $d_\ell = 10$ single dimensions per payer $\ell=1,\ldots,L$, a regularization parameter for the interpolation conditions of $0.001$, a learning rate of $\lambda=0.01$, a batch size of $100$ and $500$ training epochs. Moreover, we increase the maximum number of greedy centers to $1000$, since the dataset size increases by a factor of $500$.

For the DT ReLU NNs, we perform $1000$ training epochs with a batch size of $12$ and choose during the cross-validation between a learning rate of $\lambda=0.01,0.005,0.001,0.0005$ or $0.0001$, weight decay with parameter $\omega=0.0001$ and dropout with a dropout rate of $0.2$ as boolean hyperparameters. Further, we construct the NN architecture by choosing an increasing sequence of hidden neurons to finally match the output size of $500$ in the last linear layer. 
Again, we do not perform a cross validation for the CT ReLU NNs but empirically choose a suitable architecture using $20$ hidden neurons per hidden layer, a learning rate of $\lambda = 0.005$, a batch size of $100$, a L$2$-regularization parameter of $0.0001$ and $1000$ epochs.
For more details and the final architectures we refer to the accompanying repository.

\Cref{fig:BT_Curve} presents the relative test error comparison for all three train-test data splits between the fine-tuned deep VKOGA models and the ReLU NNs with respect to an increasing number of layers. 
As can be seen, the deep VKOGA models achieve always smaller relative test errors than the ReLU NNs, with the sole exception being $8$L-DT-VKOGA for train-test split $1$. 
Even the $1$L-VKOGA models achieve more accurate predictions compared to all NN models.
For train-test split $1$, the most accurate model is the $2$L-DT-VKOGA model, whereas the $4$L-DT-VKOGA model is the most accurate model for split $2$ and the $8$L-DT-VKOGA model for split $3$, demonstrating again the conceptual advantage of deeper kernel architectures.
Compared to the DT-VKOGA models, the CT-VKOGA models often achieve similar accuracies, but appear more stable when changing the number of layers in the sense that the relative test error usually remains in the same region of magnitude, excluding the single layer models.
Moreover, we observe that the relative test errors for train-test split $3$ are in general larger for most models compared to splits $1$ and $2$. 
This behavior is caused by a single outlier breakthrough curve, which is only contained in the test dataset for train-test split $3$, impacting the resulting relative test error significantly.

Examining the corresponding offline and online costs in \cref{tab:BTcurve_test}, we can draw similar conclusions compared to \cref{sec:5.1}, at least in the discrete-time case: 
The deep DT-VKOGA models achieve always slightly smaller offline and online runtime costs when compared to the ReLU NNs with the same number of optimizable layers. 
The single exception is the $6$L-DT-VKOGA model, which is slightly less efficient in the online phase compared to the $6$L-DT-NN.
Examining the continuous-time case, we see a different picture. First, the CT models require significantly more computational resources in the offline phase, simply due to the increased training dataset sizes and, in the case of the deep CT-VKOGA models, the increased number of greedy iterations. In detail, a single training epoch requires significantly more batches for larger datasets and every greedy iteration is also more expensive due to the need to compute the greedy selection criterion at each iteration for every training sample. Second, while the offline runtime is still comparable between the CT-NNs and the CT-VKOGA models, the online runtime of the CT-VKOGA models is much larger compared to their NN counterparts. In fact, the CT-NNs achieve a similarly efficient online procedure compared to the DT-NNs, whereas the online costs of the CT-VKOGA models are in most cases worse by a factor of $100$ compared to the DT-VKOGA models. The reason is, that the final deep kernel evaluation, which requires distance computations between all test samples and greedy centers, becomes increasingly expensive with increasing test dataset size, especially if more greedy centers are required and if the propagating function transforms the input into a higher-dimensional feature space.

\begin{figure}[!htbp]
    \centering
    \input{BTCurve_error_rel_L.pgf}
    \caption{Mean relative test error evaluated on the fixed test dataset as a function of network depth (number of layers). The error bars indicate the $5$th-$95$th percentile range of relative test errors across the test dataset. The dashed lines correspond to the deep VKOGA models, the solid lines to the ReLU NNs. Left: Train-test split $1$. Middle: Train-test split $2$. Right: Train-test split $3$.}
    \label{fig:BT_Curve}
\end{figure}

\begin{table}[!htbp]
\centering
		\begin{tabular}{l?cc?cc?cc?}

            & \multicolumn{2}{c?}{\cellcolor{light-gray} Split 1} & \multicolumn{2}{c?}{\cellcolor{light-gray} Split 2} & \multicolumn{2}{c?}{\cellcolor{light-gray} Split 3} \\

			& \cellcolor{lighter-gray} $t_\text{off}$ in [s] & \cellcolor{lighter-gray} $t_\text{on}$ in [s] & \cellcolor{lighter-gray} $t_\text{off}$ in [s] & \cellcolor{lighter-gray} $t_\text{on}$ in [s] & \cellcolor{lighter-gray} $t_\text{off}$ in [s] & \cellcolor{lighter-gray} $t_\text{on}$ in [s]\\
			\hline
            \hline
            \cellcolor{light-gray}1L-DT-VKOGA & 5.70e-03 & 1.50e-04 & 6.33e-03 & 2.14e-05 & 5.79e-03 & 2.87e-05 \\ 

            \cellcolor{light-gray}2L-DT-VKOGA & 6.38e-01 & 5.76e-04 & 8.96e-01 & 5.25e-05 & 6.55e-01 & 6.24e-05 \\ 

            \cellcolor{light-gray}4L-DT-VKOGA & 1.17e+00 & 1.06e-04 & 1.69e+00 & 2.11e-04 & 1.45e+00 & 3.94e-04 \\ 

            \cellcolor{light-gray}6L-DT-VKOGA & 1.71e+00 & 2.38e-04 & 2.55e+00 & 1.32e-04 & 2.26e+00 & 9.04e-05 \\ 

            \cellcolor{light-gray}8L-DT-VKOGA & 2.91e+00 & 2.59e-04 & 2.97e+00 & 2.63e-04 & 2.90e+00 & 1.50e-04 \\ 
            \hline
            \hline
            \cellcolor{light-gray}1L-CT-VKOGA & 3.30e+01 & 1.65e-02 & 3.24e+01 & 1.72e-02 & 3.24e+01 & 1.71e-02 \\ 

            \cellcolor{light-gray}2L-CT-VKOGA & 1.22e+02 & 2.38e-02 & 1.47e+02 & 1.82e-02 & 1.16e+02 & 1.86e-02 \\ 

            \cellcolor{light-gray}4L-CT-VKOGA & 2.37e+02 & 2.35e-02 & 2.72e+02 & 2.17e-02 & 2.89e+02 & 2.15e-02 \\ 

            \cellcolor{light-gray}6L-CT-VKOGA & 4.17e+02 & 2.48e-02 & 3.99e+02 & 2.56e-02 & 4.33e+02 & 2.33e-02 \\ 

            \cellcolor{light-gray}8L-CT-VKOGA & 5.58e+02 & 2.62e-02 & 6.00e+02 & 2.54e-02 & 5.98e+02 & 2.56e-02 \\ 
            \hline
            \hline
            \cellcolor{lighter-gray}2L-DT-NN & 1.85e+00 & 3.34e-04 & 1.34e+00 & 1.82e-04 & 1.38e+00 & 1.67e-04 \\ 

            \cellcolor{lighter-gray}4L-DT-NN & 3.25e+00 & 4.26e-04 & 2.84e+00 & 3.49e-04 & 2.79e+00 & 3.04e-04 \\ 

            \cellcolor{lighter-gray}6L-DT-NN & 2.77e+00 & 1.75e-04 & 5.12e+00 & 4.92e-04 & 4.77e+00 & 4.13e-04 \\ 

            \cellcolor{lighter-gray}8L-DT-NN & 6.17e+00 & 7.25e-04 & 6.05e+00 & 5.83e-04 & 5.96e+00 & 6.27e-04 \\ 
            \hline
            \hline
            \cellcolor{lighter-gray}2L-CT-NN & 6.87e+01 & 5.02e-04 & 9.64e+01 & 4.60e-04 & 6.84e+01 & 5.17e-04 \\ 

            \cellcolor{lighter-gray}4L-CT-NN & 1.04e+02 & 1.46e-04 & 1.26e+02 & 1.20e-04 & 1.05e+02 & 9.96e-05 \\ 

            \cellcolor{lighter-gray}6L-CT-NN & 1.36e+02 & 1.29e-04 & 1.38e+02 & 1.46e-04 & 1.37e+02 & 9.86e-05 \\ 

            \cellcolor{lighter-gray}8L-CT-NN & 1.73e+02 & 1.69e-04 & 1.72e+02 & 1.67e-04 & 1.73e+02 & 1.59e-04 \\ 
            \hline
            \hline
		\end{tabular}
	\caption{Averaged online and offline runtime measurements in seconds for all deep kernel models and NNs applied to the breakthrough curve approximation task with train-test splits $1$,$2$ and $3$.}
	\label{tab:BTcurve_test}
\end{table}

Finally, we compare the best performing deep VKOGA models and NNs on all three train-test splits to novel Graph Neural Networks (GNNs). 
Our motivation is to compare the deep VKOGA models not only to commonly used NNs but also to the new and upcoming method of spatio-temporal GNNs for the solution of spatio-temporal problems \cite{Sahili2023}. 
In the context of the breakthrough curve approximation problem we identify the parameter domain as spatial component. 
The underlying architecture of the graph is based on the combination of the parameter domain with the temporal component by edges inside the graph \cite{Kapoor2020, Roy2021}. 
On this architecture we train a graph attention network (GAT) \cite{Veli2018}  as a special case of a GNN, which is characterized by learnable attention weights between nodes. 
The applied GNN belongs to the class of CT surrogate models as defined in \cref{eq:CT_surrogate}. 
In detail, we perform $1000$ training epochs with a batch size of $1$ using a GNN training dataset of $20$ randomly created graphs from the original training data. 
To setup the GNN we start with four dense neural network layers connected to three GNN layers. The GNN layers are itself connected by 3 dense neural network layers and finally we have again four dense neural network layers to receive the output of the whole GNN. In the dense neural network layers we switch between $128$ and $256$ hidden neurons and use a learning rate of $\lambda=0.01$ and a learning rate decay of $0.8$ after every $60$-th epoch.
For more details about the used GNN we refer to the DaRUS dataset.

The relative test errors for all models on all three splits are visualized in \cref{fig:BT_Curve_pODE}. 
As expected, we see that the deep VKOGA models outperform all other models in terms of accuracy on all three train-test splits.
While the GNNs are also less accurate compared to the deep VKOGA models, they outperform the CT-NNs on all three splits and the DT-NN on split $1$, at least in terms of accuracy. 
However, in terms of offline and online efficiency, the GNNs are by far the most inefficient models in our comparison with online runtimes around $3.88e+00$ seconds and offline runtimes larger than $2.8e+03$ seconds on all three splits.

\begin{figure}[!htbp]
    \centering
    \input{BTCurve_pODE_error_rel.pgf}
    \caption{Mean relative test error evaluated on the test dataset as a function of the train-test split. The error bars indicate the $5$th-$95$th percentile range of relative test errors across the test dataset. The dashed lines correspond to the deep VKOGA models, the solid lines to the ReLU NNs and GNNs.}
    \label{fig:BT_Curve_pODE}
\end{figure}

Summarizing the numerical experiments regarding the breakthrough curve approximation problem, the DT- and CT-VKOGA models achieve in general smaller relative test errors compared to their NN counterparts as well as compared to spatio-temporal GNNs. Further, the DT-VKOGA models are even the most efficient models in both the offline and online phase.

\subsection{Parameterized ODEs}\label{sec:5.3}

\paragraph{The parameterized ODE approximation problem}
The last problem class we examine is concerned with the approximation of parameterized ODE solutions.
In particular, we employ the deep VKOGA models, ReLU NNs and GNNs onto datasets generated for the Lotka-Volterra ODE and the Brusselator ODE. 

The Lotka-Volterra ODE describes a dynamical system of predator and prey populations by the first order nonlinear ODE system 
\begin{align*}
    \frac{d}{dt} u^{(1)}(t;\mu) &= \mu_1 u^{(1)}(t;\mu) - u^{(1)}(t;\mu) u^{(2)}(t;\mu) \\
    \frac{d}{dt} u^{(2)}(t;\mu) &= - \mu_2 u^{(2)}(t;\mu) + u^{(1)}(t;\mu) u^{(2)}(t;\mu) .
\end{align*}
Here, $u^{(1)}$ and $u^{(2)}$ describe the evolution of the populations corresponding to the prey and 
predator species, respectively, over the time interval $\Omega_t \coloneqq [0,30]$, where we prescribe the parameter-independent initial conditions 
as $u^{(1)}(t=0;\mu) = 2$ and $u^{(2)}(t=0;\mu) = 4$. 
The parameters $\mu = (\mu_1,\mu_2)^T \in \Omega_\mu \coloneqq [0.8,1.2]^{2}$ describe the 
population growth and decay rates. 

The Brusselator ODE describes, e.g., symmetry breaking instable dissipative systems and autocatalytic chemical reactions. 
Simplified, the Brusselator system can be formulated as a kinetic equation system for two intermediate chemical components 
$u^{(1)}$ and $u^{(2)}$:
\begin{align*}
    \frac{d}{dt} u^{(1)}(t;\mu) &= 1 + u^{(1)}(t;\mu) \left( u^{(1)}(t;\mu) \cdot u^{(2)}(t;\mu) - 1 - \mu \right) \\
    \frac{d}{dt} u^{(2)}(t;\mu) &= u^{(1)}(t;\mu) \left( \mu - u^{(1)}(t;\mu) \cdot u^{(2)}(t;\mu) \right)
\end{align*}
We consider the time interval 
$\Omega_t \coloneqq [0,30]$, starting with $u^{(1)}(0;\mu) = u^{(2)}(0,\mu) = 1$ 
for all $\mu \in \Omega_\mu \coloneqq [0,5]$. 

\paragraph{Surrogate models}
We again distinguish between the CT and DT approach.
In the case of the DT approach, we train two surrogate models to learn the two mappings
\begin{align*}
    g^{(j)}: \mathbb{R}^P \to \mathbb{R}^{N_t}, \quad \mu \mapsto u^{(j)}(\mu) \coloneqq 
    \left[u^{(j)}(t_0,\mu),u^{(j)}(t_1,\mu),\ldots,u^{(j)}(t_{N_t-1},\mu)\right]^T
\end{align*}
for $j=1,2$ with $P$ referring to the dimension of the parameter $\mu$.
Alternatively, in the case of the CT approach, we train two surrogate models to learn the two mappings
\begin{align*}
    g^{(j)}: \Omega_t \times \mathbb{R}^P \to \mathbb{R}, \quad (t,\mu) \mapsto u^{(j)}(t;\mu), \quad j=1,2 .
\end{align*}
Note, that the distinction between a separate surrogate model for each equation of the ODE system allows for more fine-tuned architecture selections and, at least in our numerical experiments, for slightly more accurate predictions. 
In addition, GNNs correspond to the class of CT surrogate models, but do not distinguish between the different equations of the ODE system, that is a single GNN predicts the solution of both equations simultaneously.

\paragraph{Numerical experiments}
To assemble three different training datasets, we choose $N_\text{train} \in \{36,64,100\}$ equidistantly distributed points in the parameter domain $\Omega_\mu$, 
whereas the test dataset always consists of $100$ random uniformly distributed points in $\Omega_\mu$. 
To create the corresponding data labels, we consider the equidistant grid in time 
\begin{align*}
    \Omega_t^{\text{h}} \coloneqq \left\{ t_i = i \cdot 0.1,\quad i = 1,2,\ldots,N_t^\text{h}-1 \right\}, \quad 
    N_t^\text{h} \coloneqq 301
\end{align*}
and compute the numerical solutions of the ODEs using SciPy's adaptive timestepping \textit{odeint} and linear 
interpolation to infer to the values on the given grids.
Further, we adjust the relative test error measure to account for the two equations in the ODE systems as follows:
\begin{align*}
    e^\text{rel} \coloneqq \frac{1}{2 N_\text{test}} \sum_{n=1}^{N_\text{test}} \sum_{j=1}^{2} 
    \frac{\Vert u^{(j)}(\mu^{(n)}) - s(\mu^{(n)}) \Vert_2^2}{\Vert u^{(j)}(\mu^{(n)}) \Vert_2^2}
\end{align*}
Here, $u^{(j)}(\mu^{(n)}) \in \mathbb{R}^{N_t^\text{h}}$ denotes the solution of the $j$-th ODE equation for parameter sample $\mu^{(n)}$, $1 \leq n \leq N_\text{train}$ evaluated on $\Omega_t^{\text{h}}$ and $s: \Omega_\mu \to \mathbb{R}^{N_t^\text{h}}$ refers to the surrogate model.

For the ReLU NNs, we perform $1000$ training epochs with a batch size of $32$. 
In the cross-validation, we choose between $20$, $50$ or $100$ hidden neurons per hidden layer, $L=3,6$ or $9$ layers, a learning rate of $\lambda=0.01$ or $0.001$ and the use of dropout with a dropout rate of $0.2$.

For the deep VKOGA models, we choose every $4$-th training sample for DT-VKOGA and every $N_t^h$-th trainings sample for CT-VKOGA as inner first-layer centers and perform a maximum number of $1000$ greedy iterations and $1000$ training epochs with a batch size of $100$ and a regularization parameter of $\gamma = 0.0001$ for the computation of Rippa's loss.
Further, we allow the greedy algorithm to terminate early whenever the norm of the training residual falls below a tolerance of $10^{-12}$.
During the cross-validation, we choose between $d_\ell = 20$ and $50$ single dimensions per layer $\ell=1,\ldots,L$, between $L=1,2$ or $4$ layers and between a learning rate of $\lambda=0.01$ or $0.001$. 

In the case of the CT-VKOGA models and CT-NNs, we only perform a cross-validation on the smallest datasets and use the same hyperparameters on the other datasets due to the increased computational costs. 
For the same reason, we do not perform a cross-validation for the GNNs but reuse the architecture from the breakthrough curve approximation problem, see \cref{sec:5.2}. 
All model architectures are again documented in the accompanying repository.

Finally, the relative test error comparison between the DT and CT deep VKOGA models, the DT and CT ReLU NNs and the GNNs is visualized in \cref{fig:PODE_error}. 
In contrast to the relative test error comparisons from \cref{fig:MD_test,fig:BT_Curve}, we no longer examine different architecture depths but rather different training dataset sizes. 
The architecture depths of the deep VKOGA models and NNs are adjusted during the cross-validations.

\begin{figure}[!htbp]
    \centering
    \input{PODE_error_rel.pgf}
    \caption{Mean relative test error evaluated on the fixed test dataset as a function of network training dataset size $N_\text{train}$. The error bars indicate the $5$th-$95$th percentile range of relative test errors across the test dataset. The dashed lines correspond to the deep VKOGA models, the solid lines to the ReLU NNs and the GNN. Left: Lotka-Volterra ODE. Right: Brusselator ODE.}
    \label{fig:PODE_error}
\end{figure}

Starting with the DT models, we observe in most cases a slight advantage of the DT-VKOGA models compared to the DT-NNs in terms of accuracy, especially on the larger training datasets. 
Here, a deep kernel with $L=2$ layers has been chosen as a result of the cross-validation on all three training datasets for the Lotka-Volterra ODE, whereas a deep kernel with $L=4$ layers is the most frequent choice for the Brusselator ODE. 
Conversely, the fine-tuned ReLU NNs consist in most cases of $100$ hidden neurons per hidden layer and $L=6$ or $L=9$ layers.

With respect to the CT surrogate models, we observe that the GNNs yield more accurate predictions compared to the CT-NNs, which use in most cases $L=9$ layers with $50$ hidden neurons per layer.
However, the most accurate surrogate models for both parameterized ODE systems are by far the CT-VKOGA models, surpassing the accuracy of the DT surrogates and GNNs generally by more than an order of magnitude. 
In general, CT-VKOGA models with $L=4$ layers and $20$ hidden neurons per layer are selected during the cross-validation.
Exemplarily, \cref{fig:PODE_LV_NC36_CTVKOGA,fig:PODE_B_NC36_CTVKOGA} show numerical solutions and predictions of the CT-VKOGA surrogate model applied to the Lotka-Volterra and Brusselator ODE, respectively, using $N_\text{train}=36$ training data samples. 
Here, the left plot visualizes the numerical solution and prediction corresponding to the best case prediction, that is the prediction associated with the smallest relative error on the test dataset. 
Similarly, the right plot visualizes the numerical solution and prediction corresponding to the worst case prediction, whereas the middle plot corresponds to the prediction associated with the median relative test error, which separates the set of sorted relative test errors into two halves.
As the plots clearly demonstrate, the CT-VKOGA surrogate model is even on the smallest training datasets able to capture the numerical solutions remarkably well. 
Only in the worst case scenarios, small deviations between numerical solutions and predictions can be observed. 
For prediction examples corresponding to the other surrogate models and datasets, we refer to the accompanying repository.

\begin{figure}[!htbp]
    \centering
    \input{PODE_LotkaVolterraODE_NC36_CT_VKOGA.pgf}
    \caption{Ground truth and prediction of the CT-VKOGA model applied to the Lotka-Volterra ODE using $N_\text{train}=36$ training data samples and $N_\text{test}=100$ test data samples. The index $j=1,2$ refers to the solution component $u^{(j)}$ of the ODE system. Left: Best case prediction (smallest relative test error). Middle: Median case prediction (median relative test error). Right: Worst case prediction (largest relative test error).}
    \label{fig:PODE_LV_NC36_CTVKOGA}
\end{figure}

\begin{figure}[!htbp]
    \centering
    \input{PODE_BrusselatorODE_NC36_CT_VKOGA.pgf}
    \caption{Ground truth and prediction of the CT-VKOGA model applied to the Brusselator ODE using $N_\text{train}=36$ training data samples and $N_\text{test}=100$ test data samples. The index $j=1,2$ refers to the solution component $u^{(j)}$ of the ODE system. Left: Best case prediction (smallest relative test error). Middle: Median case prediction (median relative test error). Right: Worst case prediction (largest relative test error).}
    \label{fig:PODE_B_NC36_CTVKOGA}
\end{figure}

Examining the computational efficiency of all involved surrogate models, 
\cref{tab:PODE_runtime} displays the corresponding online and offline runtime measurements. There, the offline runtime measurements of the GNNs are estimations only, since the GNNs are in both offline and online phases by far the most expensive surrogate models. 
In contrast, both DT approaches yield the most efficient surrogate models in the offline and online phases, where the DT-VKOGA models are in general slightly more efficient in the offline phases but slightly less efficient in the online phases. 
Inspecting the CT surrogates, we clearly observe that the CT-VKOGA models are usually less efficient than the CT-NNs in both, offline and online phases, since the greedy iterations become increasingly expensive on larger datasets. Note, that one could still fine-tune the CT-VKOGA models, for example by decreasing the number of greedy iterations. Consequently, we expect less accurate predictions but more efficient training and inference procedures.

\begin{table}[!htbp]
  \centering
		\begin{tabular}{ll?cc?cc?}

            & & \multicolumn{2}{c?}{\cellcolor{light-gray} Lotka-Volterra ODE} & \multicolumn{2}{c?}{\cellcolor{light-gray} Brusselator ODE} \\

			& & \cellcolor{lighter-gray} $t_\text{off}$ in [s] & \cellcolor{lighter-gray} $t_\text{on}$ in [s] & \cellcolor{lighter-gray} $t_\text{off}$ in [s] & \cellcolor{lighter-gray} $t_\text{on}$ in [s] \\
			\hline
            \hline
        \multirow{5}{*}{$N_\text{train}=36$}& \cellcolor{light-gray} DT-VKOGA & 1.16e+00 & 7.77e-03 & 1.02e+00 & 1.36e-04 \\ 
        & \cellcolor{light-gray} CT-VKOGA & 5.64e+02 & 7.16e-02 & 4.79e+02 & 4.79e-02 \\ 
        & \cellcolor{lighter-gray} DT-NN & 2.32e+00 & 5.82e-04 & 1.14e+00 & 6.05e-05 \\ 
        & \cellcolor{lighter-gray} CT-NN & 1.26e+02 & 9.23e-04 & 1.26e+02 & 8.80e-04 \\ 
        & \cellcolor{lighter-gray} GNN & $>$4.00e+04 & 1.07e+01 & $>$4.00e+04 & 9.36e+00 \\ 
        \hline
        \hline
        \multirow{5}{*}{$N_\text{train}=64$}& \cellcolor{light-gray} DT-VKOGA & 1.52e+00 & 4.16e-05 & 3.58e+00 & 5.92e-05 \\ 
        & \cellcolor{light-gray} CT-VKOGA & 8.04e+02 & 6.30e-02 & 1.11e+03 & 4.03e-02 \\ 
        & \cellcolor{lighter-gray} DT-NN & 2.80e+00 & 4.07e-05 & 2.51e+00 & 4.70e-05 \\ 
        & \cellcolor{lighter-gray} CT-NN & 2.85e+02 & 1.63e-04 & 2.02e+02 & 1.60e-04 \\ 
        & \cellcolor{lighter-gray} GNN & $>$4.00e+04 & 1.08e+01 & $>$4.00e+04 & 9.55e+00 \\ 
        \hline
        \hline
        \multirow{5}{*}{$N_\text{train}=100$}& \cellcolor{light-gray} DT-VKOGA & 2.32e+00 & 3.41e-05 & 3.14e+00 & 9.66e-04 \\ 
        & \cellcolor{light-gray} CT-VKOGA & 2.85e+03 & 4.64e-02 & 2.97e+03 & 4.33e-02 \\ 
        & \cellcolor{lighter-gray} DT-NN & 3.47e+00 & 1.12e-04 & 2.67e+00 & 3.77e-05 \\ 
        & \cellcolor{lighter-gray} CT-NN & 3.66e+02 & 1.63e-04 & 3.59e+02 & 7.83e-05 \\ 
        & \cellcolor{lighter-gray} GNN & $>$4.00e+04 & 1.13e+01 & $>$4.00e+04 & 9.53e+00 \\ 
        \hline
        \hline
		\end{tabular}
	\caption{Averaged online and offline runtime measurements in seconds for all deep kernel models and NNs applied to the parameterized ODE approximation task using training datasets of sizes $N_\text{train}=36,64,100$.}
	\label{tab:PODE_runtime}
\end{table}

\section{Conclusion and Outlook}\label{sec:6}

We underline the two main contributions of this work. 
First, we have extended the work on $2$L-VKOGA \cite{wenzel2024application,wenzel2024data,herkert2024greedy} by introducing and examining a combination of VKOGA with even deeper kernels, in particular with up to $8$ layers, which provides a more powerful class of function approximators. The resulting deep VKOGA models inherit not only the sparsity, stability and theoretical framework of greedy kernel interpolation, but also the flexibility, hierarchical feature learning and excellent approximation capabilities of classical NNs, rendering them widely and efficiently applicable. Second, we have systematically compared the deep VKOGA models to fully-connected ReLU NNs with a focus on prediction accuracy, offline and online efficiency and network depth. In particular, three representative problem classes were considered: The approximations of model problems with increasingly higher input dimensions and increasingly difficult characteristics, the approximations of breakthrough curves of chemical species flowing through 3D porous geometries and the approximations of parameterized ODE solutions. To approximate parameterized ODE solutions, we have distinguished between a continuous- and discrete-time approach and also compared the deep VKOGA models to another class of NNs, namely GNNs, which have specifically been adjusted for the task of time-series approximations. Across all settings and problem classes, the deep VKOGA models have proven to consistently outperform the NNs in terms of accuracy and, in many cases, also in terms of computational efficiency, highlighting the potential of deep VKOGA models for a large variety of applications. Further, we have showcased that deeper kernel architectures can be advantageous, especially for complex target functions with different characteristics.

Despite the imminent success of the deep VKOGA models, our conducted experiments have also revealed a major limitation: For very large datasets, the deep kernel training, the greedy iterations and the evaluation in the online phase require an increased amount of computational resources compared to all tested NNs. The reasons are that Rippa's loss is much more expensive than the MSE loss, that the training residual must be evaluated on the complete training dataset to select the new greedy center at each greedy iteration and that the computation of the deep kernel evaluations is more expensive for a larger number of selected greedy centers. Thus, limiting the number of greedy centers, filtering the datasets a priori or changing the approach, for example from continuous-time to discrete-time predictions, should be under consideration if the efficiency needs to be improved.

Moreover, several improvements and future research topics have emerged from our study. 
First, our experiments have been limited to a small amount of architecture- and hyperparameter-combinations tested in the cross-validations. 
As a consequence, there are certainly different architectures leading to more accurate and more efficient predictions for both, NNs and deep VKOGA models. 
Second, the deep kernel is trained using a predefined number of randomly chosen fixed inner first-layer centers. 
However, there are other possible strategies to select these centers, where two are still under investigation: 
On the one hand, the fixed inner first-layer centers could coincide with the greedy centers, resulting in a deep kernel that needs to be updated iteratively as well. 
We expect that this increases the computational costs of the offline phase of the corresponding deep VKOGA model, but also provides interpretability in the scope of the deep kernel representer theorem \cite{bohn2019representer}. 
On the other hand, the fixed inner first-layer centers could be given as optimizable parameters to the deep kernel, which are adjusted during training. 
This leads to a flexible and problem-adapted choice, but requires careful regulations to ensure that the inner centers remain meaningful and do not collapse.
Third, the combination of convolutional NNs and deep kernels leading to convolutional deep kernels represents another promising avenue for future research. 
In particular when processing tensor-valued input data, convolutional kernel layers might be able to process the spatial information directly. 
This is particularly interesting in the framework of the breakthrough curve approximation problem, where convolutional deep kernels might replace the PCA feature map and directly handle the 3D geometry input data.
\section*{Supplementary information}

\paragraph{Acknowledgements}
We acknowledge the support by the Stuttgart Center for Simulation Science (SimTech). We thank T. Wenzel for discussions on deep kernel approaches.

\paragraph{Funding}
Funded by Deutsche Forschungsgemeinschaft (DFG, German Research Foundation) in Project No. 314733389, 540080351 and in SPP 2311 in Project No. 465243391, and under Germany’s Excellence Strategy - EXC 2075
– 390740016. 

\paragraph{Data availability}
The datasets considered in this study are available via DaRUS (\url{https://doi.org/10.18419/DARUS-5167}). 

\paragraph{Code availability}
The source code for the numerical experiments presented in this study is available via DaRUS (\url{https://doi.org/10.18419/DARUS-5167}). 

\bibliographystyle{plain}
\bibliography{template}  






\end{document}